\documentclass[10pt]{article}

\usepackage[latin1]{inputenc}
\usepackage[english]{babel}
\usepackage{amssymb}
\usepackage{amsmath}
\usepackage{enumerate}
\usepackage{indentfirst}
\usepackage{color}
\usepackage{textcomp}
\usepackage{fancyhdr}
\usepackage{titlesec}

\newtheorem{theorem}{Theorem}[section]
\newtheorem{definition}{Definition}[section]
\newtheorem{corollary}{Corollary}[section]
\newtheorem{proposition}{Proposition}[section]
\newtheorem{lemma}{Lemma}[section]

\newtheorem{theoreme}{Théorème}[section]
\newtheorem{remark}{Remark}[section]

\newcommand{\eq}[1]{\begin{equation}#1\end{equation}}
\newcommand{\arr}[2]{\begin{array}{#1}#2\end{array}}
\newcommand{\theo}[1]{\begin{theorem}#1\end{theorem}}

\newcommand{\R}{\mathbb{R}}
\newcommand{\dv}{\hbox{div} \hspace{0.1cm}}
\newcommand{\curl}{\hbox{curl} \hspace{0.1cm}}

\newcommand{\lem}[1]{\begin{lemma}#1\end{lemma}}
\newcommand{\aligne}[1]{\begin{align*}#1\end{align*}}
\newcommand{\cent}[1]{\begin{center}#1\end{center}}

\makeatletter

\@addtoreset{equation}{section}
\makeatother

\title{Asymptotic profiles for the second grade fluids equations in $\R^3$}
\author{Olivier Coulaud}
\date{}
\begin{document}
\maketitle
\begin{abstract}
In the present paper, we study the long time behaviour of the solutions of the second grade fluids equations in $\R^3$. Using scaling variables and energy estimates in weighted Sobolev spaces, we describe the first order asymptotic profiles of these solutions. In particular, we show that the solutions of the second grade fluids equations converge to self-similar solutions of the heat equations, which are explicit and depend on the initial data. Since this phenomenon occurs also for the Navier-Stokes equations, it shows that the fluids of second grade behave asymptotically like Newtonian fluids.
\end{abstract}
\section{Introduction}
Since one can find a lot of non-Newtonian fluids in the nature or in the industry, their mathematical study is a significant topic of research. For instance, wet sand or the paste used to make paper are non-Newtonian fluids. In this paper, we investigate the long time behaviour of a particular class of non-Newtonian fluids, namely the second grade fluids. The equations which describe such fluids have been introduced from a mathematical point of view in 1974 by Dunn and Fosdick in \cite{dunnfosdick74} and have been the topic of many research works in mathematics. These fluids are a particular case of a large class of non-Newtonian fluids, called fluids of differential type, or Rivlin-Ericksen fluids (see \cite{rivlinericksen55}). The constitutive laws of the differential fluids are given through the Rivlin-Ericksen tensors, defined recursively by
\cent{$
\arr{l}{
A_1 = \nabla u + \left(\nabla u\right)^t,\\
A_k = \partial_t A_{k-1} + u.\nabla A_{k-1} + \left(\nabla u\right)^t A + A \nabla u,
}
$}
where $u:\R^d \rightarrow \R^d$ is a vector field which represents the velocity of a fluid filling a domain of $\R^d$, $d =2,3$. According to this model, the equations of the fluids of grade $n \in \mathbb N$ are obtained by considering the stress tensor
\cent{
$
\sigma = - p Id + Q \left( A_1, A_2,..., A_n\right),
$
}
where $p$ is the pressure of the fluid and $Q$ is a polynomial function of degree $n$. Notice that the fluids of grade 1 correspond to the classical Navier-Stokes equations, which concern Newtonian fluids. According to the model of Dunn and Fosdick (see \cite{dunnfosdick74}), the constitutive law of the second grade fluids is obtained via the stress tensor
\cent{
$
\sigma = - p Id + \nu A_1 +\alpha_1 A_2 + \alpha_2 A^2_1,
$
}
where $\nu>0$ is the dynamic viscosity of the fluid, $\alpha_1 >0$ and $\alpha_2 \in \R$. In \cite{dunnfosdick74}, thermodynamic considerations led the authors to assume that $\alpha_2 =-\alpha_1$. Consequently, we replace $\alpha_1$ by $\alpha$. Introduced in the equations of conservation of momentum, the tensor $\sigma$ leads to the system of equations
\eq{\label{g2u}
\arr{l}{\partial_t\left(u-\alpha \Delta u\right)-\nu \Delta u +\curl\left(u-\alpha \Delta u\right)\wedge u + \nabla p=0,\\
\dv u=0,\\
u_{\left|t=0\right.}=u_0,
}
}
where $\wedge$ denotes the classical vectorial product on $\R^3$, $p$ is the pressure which depends on $u$ and $u_0$ is the initial data. In the two-dimensional case, we have used the convention that $u=\left(u_1,u_2,0\right)$ and $\curl u = \left(0,0,\partial_1 u_2-\partial_2 u_1\right)$.
\vspace{0.5cm}\\
\indent Several existence and uniqueness results have been obtained for this system of equations, mainly on a bounded domain $\Omega$ of $\R^2$ or $\R^3$ with Dirichlet or periodic boundary conditions (see for instance \cite{bernard99}, \cite{breschlemoine98}, \cite{cioranescuelhacene84}, \cite{cosciagaldi93}, \cite{cosciagaldi94}, \cite{cosciagaldi00}, \cite{galdisequeiravideman97}, \cite{novotnysequeiravideman97} or \cite{leroux99}). The first existence and uniqueness result has been obtained by Cioranescu and El Hacène in 1984 in \cite{cioranescuelhacene84}. They have shown, on a bounded set of $\R^d$, $d=2,3$, with homogeneous boundary conditions, that there exists a unique weak solution to (\ref{g2u}) belonging to the space $L^\infty\left(\left[0,T\right],H^3(\Omega)^d\right)$, where $T>0$ and $H^s(\Omega)$ denotes the Sobolev space of order $s$ (see \cite{cioranescuelhacene84}). Besides, this solution is global in time when the space dimension is $2$. This result is based on a priori estimates and a Galerkin approximation with a basis of eigenfunctions corresponding to the scalar product associated to the operator $\curl\left(u-\alpha \Delta u\right)$.  In the same case, using the Schauder fixed point theorem, Galdi, Grobbelaar-Van Dalsen and Sauer established the existence and uniqueness of classical solutions to (\ref{g2u}) when the data belong to $H^m$, with $m\geq 5$ (see \cite{galdivandalsensauer93}). They also have shown that these solutions are global in time, provided that the initial data are small enough in $H^m(\Omega)$. Later, Cioranescu and Girault improved the results of \cite{cioranescuelhacene84} and \cite{galdivandalsensauer93} and showed that the local weak solutions belonging to $H^3(\Omega)$ are actually global in time in dimension 3 if the data are small enough and are strong solutions if the data belong to $H^m$, $m\geq 4$ (see \cite{cioranescugirault97}). Finally, Bresch and Lemoine have generalized the results of \cite{galdivandalsensauer93}, \cite{cioranescuelhacene84} and \cite{cioranescugirault97} in dimension 3 in establishing the existence and uniqueness of local solutions belonging to the space $W^{2,r}(\Omega)$ with $r>3$. Furthermore, they have shown that these solutions are global in time if the initial data are small enough in $W^{2,r}(\Omega)$ (see \cite{breschlemoine98}). In this work, instead of applying a Galerkin approximation, the authors used Schauder's fixed point Theorem.
\vspace{0.5cm}\\
\indent In the present paper, we are interested in the description of the asymptotic profiles of the solutions of second grade fluids equations. In what follows, we consider a second grade fluid which fills the whole space $\R^3$, without any forcing term applied to it. In this case, if the initial data are small enough, the solutions of such a system tend to $0$ when the time $t$ goes to infinity. The aim of this study is to investigate the way that these solutions go to $0$. More precisely, we will show that the solutions of (\ref{g2u}) behave asymptotically like self-similar solutions to the heat equation, which are smooth and that one can compute explicitly from the data. In this article, we restrict ourselves to the study of the first order asymptotic profile, that is to say that the speed of the convergence of the solutions of (\ref{g2u}) to explicit smooth functions is limited by spectral considerations. For the Navier-Stokes equations, there already exist several results that describe the asymptotic profiles of the solutions. In dimension 2 and 3, Gallay and Wayne have shown in \cite{gallaywayne02} and \cite{gallaywayne02bis} that the first order asymptotic profiles of the solutions of the Navier-Stokes equations are given up to a constant by smooth Gaussian functions which are self-similar solutions to the heat equations. These results hold with restrictions on the size of the data, but, in dimension 2, the convergence has been generalized to the case any data in \cite{gallaywayne05}. For this work, the authors applied arguments that come from the study of dynamical systems. In fact, they have shown the existence of a finite-dimensional manifold locally invariant by the semiflow associated to the Navier-Stokes equations. Then, they proved that the solutions of the Navier-Stokes equations are locally attracted by this manifold, and consequently behave like the solutions on it. The study of the dynamics of the Navier-Stokes equation onto this manifold gave them the description of the first and second order asymptotic profiles. The asymptotic profiles of the solutions of the equations of second grade fluids have been studied in $\R^2$ by Jaffal-Mourtada in \cite{jaffal11}. She has shown, under smallness assumptions on the data, that the first order asymptotic profiles of the solutions of the second grade fluids equations are the same as the ones described by Gallay and Wayne in \cite{gallaywayne02} for the Navier-Stokes equations. However, the method that she used in \cite{jaffal11} is slightly different from the one used in \cite{gallaywayne02}. Indeed, instead of showing the existence of an invariant manifold, the author performed energy estimates in various function spaces, notably weighted Sobolev spaces. The concrete interpretation of this result is that, in dimension 2, the fluids of second grade behave asymptotically like Newtonian fluids. In this article, we are interested in the generalization of this result to the dimension 3. Notice that there are significant differences in the asymptotic behaviour of the Navier-Stokes equations between the cases of $\R^2$ and $\R^3$. Indeed, in dimension 2, the asymptotic profiles of the Navier-Stokes equations are given up to a constant by a Gaussian function called the Oseen vortex sheet. In dimension 3, the first order asymptotic profiles of the solutions are defined as the linear combination of three distinct smooth functions (see Section \ref{sectheo}).
\vspace{0.5cm}\\
Actually, the system that we study in this article is not exactly (\ref{g2u}) but the one satisfied by the vorticity $w =\curl u$. The motivation to do this comes from the fact that, due to spectral considerations which will be explained more precisely later, we have to solve the equations of second grade fluids in weighted Lebesgue spaces. Unfortunately, the system (\ref{g2u}) do not preserve in general the weighted Lebesgue spaces. We assume, for the sake of simplicity, that $\nu=1$ and consider initial vorticity data $w_0$. Taking formally the $\curl$ of (\ref{g2u}), we get the vorticity system of equations
\eq{\label{g2w}
\arr{l}{\partial_t\left(w-\alpha \Delta w\right)- \Delta w + \curl \left(\left(w-\alpha \Delta w\right)\wedge u\right) =0,\\
\dv u=\dv w=0,\\
w_{\left|t=0\right.}=w_0.}
}
This system is actually autonomous. Indeed, provided that $w$ is sufficiently smooth, the divergence free vector field $u$ can be recovered from $w$ via the Biot-Savart law, which is a way to get a divergence free vector field from its given vorticity. It is defined by the formula
\eq{
\displaystyle u(x)=-\frac{1}{4\pi}\int_{\R^3} \frac{\left(x-y\right)\wedge w(y)}{\left|x-y\right|^3} dy.
}
In Section \ref{sectheo}, more details are given on the Biot-Savart law and its property (see Lemma \ref{biots}).
In this article, we show that the solutions of the system (\ref{g2w}) behave asymptotically like vector fields whose components are self-similar solutions to the well known heat equations, that is to say under the form
\cent{
$
\displaystyle\left(t,x\right) \rightarrow \frac{1}{\left(t+T\right)^2} F\left(\frac{x}{\sqrt{t+T}}\right),
$
}
where $F$ is a vector field of $\R^3$ and $T$ is a positive constant.
\vspace{0.5cm}\\
\indent We introduce now a powerful tool in the study of the asymptotics of solutions to partial differential equations, that is scaled variables or self-similar variables. In order to define those variables, we set a positive constant $T$, and we will always assume $T\geq 1$. The motivation to introduce this constant is that, by doing this, we will be able to establish the convergence of the solutions to their asymptotic profiles without restriction on the size of the constant $\alpha$. As it is explained below, the constant $T$ will be chosen large enough to have $\frac{\alpha}{T}$ small enough. For a solution $w$ of the system (\ref{g2w}), we define $W$ and $U$ through the change of variable $\displaystyle X=\frac{x}{\sqrt{t+T}}$ and $\tau = \log (t+T)$. More precisely, we set
\eq{\label{scaledvar2}\arr{l}{
\displaystyle w\left(t,x\right)= \frac{1}{t+T} W \left( \log \left(t+T\right), \frac{x}{\sqrt{t+T}} \right),\\
\displaystyle u\left(t,x\right)= \frac{1}{\sqrt{t+T}} U \left( \log \left(t+T\right), \frac{x}{\sqrt{t+T}} \right).}
}
Equivalently, we have the equalities
\eq{\label{scaledvar1}\arr{l}{
W\left(\tau,X\right)= e^\tau w \left( e^\tau - T, e^{\tau/2} X\right),\\
U\left(\tau,X\right)= e^{\tau/2} u \left( e^\tau - T, e^{\tau/2} X\right).}} 
\indent Scaling variables have been initially introduced to study the asymptotic behaviours of solutions of parabolic equations, and in particular to show the convergence to self-similar solutions (see  \cite{escobedokavianmatano95}, \cite{escobedozuazua91}, \cite{galaktionov91} or \cite{kavian87}). Actually, this tool is also efficient to study the long-time behaviour of a lot of various equations, not necessarily parabolic ones. For instance, Gallay and Raugel used them to describe the first and second order asymptotic profiles of the solutions to damped waved equations (see \cite{gallayraugel98}) and to show the stability of hyperbolic fronts (see \cite{gallayraugel00}). Self-similar variables have been also used to study the asymptotic profiles of the Navier-Stokes equations (see \cite{gallaywayne02}, \cite{gallaywayne02bis}, \cite{gallaywayne05} and \cite{gallaywayne06}) and the second grade fluids equations in dimension 2 (see \cite{jaffal11}). Assuming that $w$ is a solution of (\ref{g2w}), a short computation shows that $W$ is a solution of the system
\eq{\label{g2W}
\arr{l}{\partial_\tau\left(W-\alpha e^{-\tau} \Delta W\right) - \mathcal L (W) +\curl\left(\left(W-\alpha e^{-\tau} \Delta W\right)\wedge U\right)\\
\hspace{7cm}+\alpha e^{-\tau}\Delta W+\alpha e^{-\tau}\frac{X}{2}.\nabla \Delta W=0,\\
\dv U=\dv W=0,\\
W_{\left|\tau=\log(T)\right.}=W_0,}
}
where $\mathcal L$ is the linear differential operator defined by
\cent{$\mathcal L(W)=\Delta W +W+\frac{X}{2}.\nabla W$.}
We first emphasize that the system (\ref{g2W}) is now non-autonomous and initialised at $\tau=\log(T)$, that is the reason why we introduced the constant $T$. Indeed, this operation allows to avoid restrictions on the size of $\alpha$ by choosing $T$ large enough. We also notice that, in the first equality, several terms formally tend to $0$ when $\tau$ goes to infinity. Actually, the main theorem of this article shows that the solutions of (\ref{g2W}) converge when $\tau$ goes to infinity to particular solutions to the equality
\eq{\label{g2lim}
\partial_\tau W_\infty = \mathcal L (W_\infty).
}
More precisely, the aim of the present paper is to decompose $W$ on the spectrum of $\mathcal L$ on an appropriate space of functions and to show that the asymptotic behaviour of $W$ is dominated by the projection of $W$ onto the eigenspace corresponding to the first eigenvalue of $\mathcal L$. Additionally, this projection satisfies the equality (\ref{g2lim}). We define now the weighted Lebesgue spaces, which are suitable for the study of the spectrum of $\mathcal L$. For every $m\in \mathbb N$, one defines $L^2(m)$, given by
\cent{
$
L^2(m)= \left\{ u \in L^2(\R^3) : \left(1+\left|x\right|^2\right)^{m/2} u \in L^2(\R^3)\right\},
$
}
where $\displaystyle \left|x\right| = \left(\sum^3_{i=1} x^2_i\right)^{1/2}$.
\vspace{0.5cm}\\
By the same way, for $m\in \mathbb N$ and $n\geq 2$, we define the weighted Sobolev spaces by
\aligne{
&H^1(m)= \left\{u \in L^2(m) : \partial_i u \in L^2(m), i \in \left\{1,2,3\right\}\right\},\\
&H^n(m)= \left\{u \in L^2(m) : \partial_i u \in H^{n-1}(m), i \in \left\{1,2,3\right\}\right\}.
}
The incompressibility condition on the vector fields $W$ and $U$ makes natural to work on the spaces
\cent{$\arr{l}{
\mathbb L^2(m)=\left\{u \in L^2(m)^3 : \dv u=0\right\},\\
\\
\mathbb H^2(m) = \left\{u \in H^2(m)^3 : \dv u=0\right\},}
$}
equipped with the norms
\cent{
$\left\|u\right\|_{L^2(m)} = \left\|\left(1+\left|x\right|^2\right)^{\frac{m}{2}} u\right\|_{L^2},$} and \cent{$\left\|u\right\|_{H^2(m)} = \left(\left\|u\right\|^2_{L^2(m)}+\left\|\nabla u\right\|^2_{L^2(m)}+\left\|\nabla^2 u\right\|^2_{L^2(m)}\right)^{1/2}$.
}
In \cite{gallaywayne02bis}, Gallay and Wayne show that the spectrum of $\mathcal L$ on $\mathbb L^2(m)$ is the union of the discrete spectrum
\cent{$\sigma_d (\mathcal L)=\left\{-\frac{1}{2}\left(k+1\right), k\in \mathbb N\right\},$}
and the continuous one
\cent{$\sigma_c(\mathcal L)=\left\{\lambda \in \mathbb C : Re(\lambda) \leq \frac{1}{4}-\frac{m}{2}\right\}.$}
In order to describe the first order asymptotic profiles of the solutions of (\ref{g2W}), we need to have at least one isolated eigenvalue in the spectrum of $\mathcal L$. Looking at $\sigma_c(\mathcal L)$, we notice that one can "push" the continuous spectrum to the left by choosing $m$ large enough. For this reason, we should work at least in the weighted space $\mathbb L^2(3)$, where $-1$ is an isolated eigenvalue of $\mathcal L$. Actually, in order to be close to the optimal rate of convergence, we prefer working in $\mathbb L^2(4)$, where the discrete spectrum is $\sigma_d(\mathcal L)=\left\{-1,-\frac{3}{2}\right\}$ and the continuous one is $\sigma_c(\mathcal L)=\left\{\lambda \in \mathbb C : Re(\lambda) \leq -\frac{7}{4}\right\}$. The main aim of this article is to show that one can decompose a solution $W$ of (\ref{g2W}) into the form
\eq{\label{decomp1}
W(\tau)=\Omega(\tau)+R(\tau),
}
where $\Omega$ is an eigenfunction of $\mathcal L$ associated to the eigenvalue $-1$ and $R$ tends to $0$ faster than $\Omega$ into $\mathbb L^2(4)$ when $\tau$ goes to infinity.
\vspace{0.5cm}\\
Since the first eigenvalue smaller than $-1$ is $-\frac{3}{2}$ , the best result that one expects is
\cent{
$\displaystyle R(\tau)=\mathcal O(e^{-\frac{3\tau}{2}})$ in $\mathbb L^2(4)$.
}
Actually, the result that we obtain holds under smallness assumptions on the size of the data in $\mathbb H^2(4)$. Besides, provided that the initial data are small enough compared to the parameters of the equations, one can choose the rate of convergence as close as wanted to the optimal one.
\section{\label{sectheo}First order asymptotics and preliminary results}
Before stating the main theorem of this paper, we have describe the eigenspace of $\mathcal L$ associated to the eigenvalue $-1$. In \cite[appendix A]{gallaywayne02bis}, they show that the multiplicity of the eigenvalue $-1$ is $3$ and that a suitable basis $\left\{f_1,f_2,f_3\right\}$ of the associated eigenspace $E_{-1}$ is given by
\eq{\label{fi}
f_i = \curl (G e_i)
, \quad i=1,2,3,}
where $\displaystyle G(X) = \frac{1}{(4\pi)^{3/2}}e^{-\frac{\left|X\right|^2}{4}}$ and $\left\{e_1,e_2,e_3\right\}$ is the canonical basis of $\R^3$.
\vspace{0.5cm}\\
Through a short computation, we see that $f_i(X)=p_i(X) G(X), \quad i=1,2,3$, where
\cent{
$\displaystyle p_1(X)= \frac{1}{2}\left(\arr{c}{0\\-X_3\\X_2}\right)$, $\displaystyle p_2(X)= \frac{1}{2}\left(\arr{c}{X_3\\0\\-X_1}\right)$ and $\displaystyle p_3(X)= \frac{1}{2}\left(\arr{c}{-X_2\\X_1\\0}\right)$.
}
In particular, the vector fields $p_i$ satisfy $\quad \dv p_i=0\quad$ and $\quad\curl p_i=e_i$. Integrating by parts, we also notice that
\eq{\label{fj}
\displaystyle\int_{\R^3} p_i(X). f_j(X) dX= \int_{\R^3}\curl (p_i (X)). \left(G(X)e_j\right) dX = \left(e_i . e_j\right)\int_{\R^3} G(X) dX=\delta_{ij}.
}
Furthermore, defining $\mathcal L^*=\Delta -\frac{X}{2}.\nabla -\frac{1}{2}$ the formal adjoint of $\mathcal L$, we check easily that 
\cent{$\mathcal L^* p_i=-p_i$.}
With the basis $\left\{f_1,f_2,f_3\right\}$, the decomposition (\ref{decomp1}) can be written 
\eq{\label{decomp2}
W(\tau)=\sum^3_{i=1} \beta_i(\tau) f_i+R(\tau),
}
where $\beta_i(\tau)\in \R$.
\vspace{0.5cm}\\
As we can see in \cite{gallaywayne02bis}, $\mathbb L^2(4)= E_{-1} \oplus \mathcal W$, where
\cent{
$
\displaystyle \mathcal W =\left\{ f \in \mathbb L^2(4) : \int_{\R^3} X_i f_j(X) dX = 0, \quad i,j=1,2,3\right\}.
$
}
Consequently, one has to choose $\beta_i$ such that $\displaystyle \int_{\R^3} X_i R_j(\tau,X) dX = 0$, for $i,j \in \left\{
1,2,3\right\} $. To this end, we set
\cent{$\displaystyle \beta_i(\tau)=\int_{\R^3} p_i(X).W(\tau, X)dX$.}
In fact, assuming that $W \in \mathbb L^2(4)$ and using the divergence free property of $W$, it is easy to check that
\aligne{
\int_{\R^3} p_1(X).W(X)dX = \int_{\R^3} X_2W_3(X)dX= -\int_{\R^3} X_3W_2(X)dX,\\
\int_{\R^3} p_2(X).W(X)dX = \int_{\R^3} X_3W_1(X)dX= -\int_{\R^3} X_1W_3(X)dX,\\
\int_{\R^3} p_3(X).W(X)dX = \int_{\R^3} X_1W_2(X)dX= -\int_{\R^3} X_2W_1(X)dX,
}
and thus, using (\ref{fj}) and the decomposition (\ref{decomp2}), we can conclude that
\cent{
$
\displaystyle \int_{\R^3} X_i R_j(X) dX =0$, for all $i,j\in \left\{1,2,3\right\}$.
}
The next lemma gives more details about $\beta_i$, and shows that the projection of $W$ onto $E_{-1}$ is actually a solution of (\ref{g2lim}).
\lem{
Let $W \in C^o\left(\left[\tau_0,T\right),\mathbb H^2(4)\right)$ be a solution of (\ref{g2W}) and let 
\cent{$\displaystyle \beta_i(\tau) = \int_{\R^3} p_i(X) . W(\tau, X)dX$.}
Then, for all $\tau \in \left[\tau_0, T\right]$,
\eq{\label{beta0}
\beta_i(\tau)= b_i e^{-\tau},
}
where $\displaystyle b_i =\int_{\R^3} p_i(X).W_0(X) dX$.
}
\textbf{Proof: }The proof of this lemma is made formally, assuming that every quantity that we consider is well defined. Actually, in the remaining of this article, we will work with regularized solutions for which the next computations are rigorous. In order to get (\ref{beta0}), we only have to show that $\beta_i$ satisfies
\eq{\label{beta1}
\partial_\tau \beta_i\left(\tau\right)=-\beta_i(\tau).
}
Performing the $L^2-$scalar product of the first equality of (\ref{g2W}) with $p_i$, we obtain
\eq{
\arr{l}{\partial_\tau \beta_i(\tau) =\alpha e^{-\tau}\left(p_i,\partial_\tau \Delta W\right)_{L^2}-\alpha e^{-\tau}\left(p_i, \Delta W\right)_{L^2}+\left(p_i,\mathcal L(W)\right)_{L^2}\\
\hspace{2cm}+\left(p_i,\curl \left(\left(W-\alpha e^{-\tau} \Delta W\right)\wedge  U\right)\right)_{L^2} -\alpha e^{-\tau}\left(p_i,\Delta W+\frac{X}{2}.\nabla \Delta W\right)_{L^2}.}
}
Integrating several times by parts, it is easy to check that
\cent{
$
\alpha e^{-\tau}\left(p_i,\partial_\tau \Delta W\right)_{L^2}=\alpha e^{-\tau}\left(p_i, \Delta W\right)_{L^2}=\alpha e^{-\tau}\left(p_i,\Delta W+\frac{X}{2}.\nabla \Delta W\right)_{L^2}=0.
$
}
Thus, integrating by parts and recalling that $\curl pi = e_i$, one has
\eq{\label{beta2}
\partial_\tau \beta_i(\tau) =-\beta_i(\tau)+\int_{\R^3} e_i. \left( \left(W(X)-\alpha e^{-\tau} \Delta W(X)\right)\wedge U(X)\right)dX.
}
It remains to show that the last term of the right hand size of (\ref{beta2}) vanishes. Noticing that $W=\curl U$, an easy computation shows, for $i\in \left\{1,2,3\right\}$,
\eq{\label{gnouf}
\left(U(X)\wedge\left(W(X)-\alpha e^{-\tau} \Delta W(X)\right)\right)_i=\frac{1}{2}\partial_i\left(\left|U\right|^2\right)-U.\nabla U_i-\alpha e^{-\tau}\left(U.\partial_i \Delta U- U.\nabla \Delta U_i\right).
}
Thus, using the divergence free property of $U$ and integrating by parts, we get
\aligne{
\int_{\R^3} e_i. \left(U(X)\wedge \left(W(X)-\alpha e^{-\tau} \Delta W(X)\right)\right)dX = -\alpha e^{-\tau}\int_{\R^3} U(X).\partial_i \Delta U(X) dX.
}
Another integration by parts yields
\aligne{
\int_{\R^3} e_i. \left(U(X)\wedge \left(W(X)-\alpha e^{-\tau} \Delta W(X)\right)\right)dX = \frac{\alpha}{2} e^{-\tau}\int_{\R^3}\partial_i \left(\left|\nabla U(X)\right|^2\right) dX=0,
}
and thus we obtain (\ref{beta1}).
\begin{flushright}
$\square$
\end{flushright}
We can now state the main theorem of this paper, which shows in particular that the first order asymptotic profile of a solution $W$ in $\mathbb H^2(4)$ of (\ref{g2W}) is the same as the first order asymptotic profile obtained for the Navier-Stokes equations.
\theo{\label{theo1}
Let $\theta$ be a fixed constant such that $0<\theta<\frac{3}{2}$ and $W_0 \in \mathbb H^2(4)$. There exist two positive constants $\gamma_0 = \gamma_0(\alpha)$ and $T_0 = T_0(\alpha,\theta) \geq 1$ such that if $ T \geq T_0$ and there exists a positive constant $\gamma\leq \gamma_0$ such that
\eq{\label{cond1}
\left\|W_0\right\|^2_{L^2(4)}+ \left\|\nabla W_0\right\|^2_{L^2}+\alpha e^{-\tau_0} \left\|\Delta W_0\right\|^2_{L^2} + \alpha^2 e^{-2\tau_0}\left\|\left|X\right|^4 \Delta W_0\right\|^2_{L^2} \leq \gamma \left(\frac{3}{2}-\theta\right)^2,
}
where $\tau_0 = \log(T)$,
\vspace{0.5cm}\\
then there exist a unique solution $W \in C^0\left(\left[\tau_0,+\infty\right), \mathbb H^2(4)\right)$ to the system (\ref{g2W}) and a positive constant $C=C(\theta, \alpha, T_0)$ such that
\eq{\label{inetheo1}
\left\|\left(I-\alpha e^{-\tau}\Delta\right)\left(W(\tau)-e^{-\tau}\sum^3_{i=1} b_i f_i\right)\right\|_{L^2(4)} \leq C \gamma \left(\frac{3}{2}-\theta\right) e^{-\theta \tau},
}
where $\displaystyle b_i = \int_{\R^3} p_i(X) . W_0( X)dX$.
}
In the classical variables, the next corollary is deduced from Theorem \ref{theo1}.
\begin{corollary}
\label{cor1}
Let $\theta$ be a constant such that $0<\theta<\frac{3}{2}$,  $w_0\in \mathbb H^2(4)$ and $\displaystyle b_i = \frac{1}{T} \int_{\R^3} p_i(x). w_0(x) dx$. There exist $\gamma_0=\gamma_0(\alpha)>0$ and $T_0=T_0(\alpha,\theta)\geq 1$ such that if there exist $T \geq T_0$ and $\gamma\leq \gamma_0$ such that
\eq{\arr{l}{
T^{1/2} \left\|w_0\right\|^2_{L^2}+T^{-7/2} \left\|\left|x\right|^4 w_0\right\|^2_{L^2}+T^{3/2}\left\|\nabla w_0\right\|^2_{L^2}\\
\hspace{4cm}+\alpha T^{3/2} \left\|\Delta w_0\right\|^2_{L^2}+\alpha^2 T^{-3/2} \left\|\left|x\right|^4\Delta w_0\right\|^2_{L^2}\leq \gamma \left(\frac{3}{2}-\theta\right)^2,}
}
then there exists a unique solution $w\in C^0\left(\left[0,+\infty\right),\mathbb H^2(4)\right)$ to the system (\ref{g2w}) such that, for all $1\leq p \leq 2$, the following inequality holds
\eq{
\left\|\left(I-\alpha \Delta\right)\left(w(t)-\sum^3_{i=1}\frac{b_i}{\left(t+T\right)^2}f_i\left(\frac{x}{\sqrt{t+T}}\right)\right)\right\|_{L^p} \leq C \gamma\left(\frac{3}{2}-\theta\right)\left(t+T\right)^{-1-\theta+\frac{3}{2p}},
}
where $C=C(\theta, \alpha, T_0)$  is a positive constant. Besides, for all $1\leq p \leq +\infty$, one has
\eq{
\left\|w(t)-\sum^3_{i=1}\frac{b_i}{\left(t+T\right)^2}f_i\left(\frac{x}{\sqrt{t+T}}\right)\right\|_{L^p} \leq C \gamma\left(\frac{3}{2}-\theta\right)\left(t+T\right)^{-1-\theta+\frac{3}{2p}}.
}
Let $u$ be the divergence free vector field obtained from $w$ through the Biot-Savart law. For all $\frac{3}{2}\leq p\leq +\infty$, one has
\eq{
\left\|u(t)-\sum^3_{i=1}\frac{b_i}{\left(t+T\right)^{3/2}}v_i\left(\frac{x}{\sqrt{t+T}}\right)\right\|_{L^p} \leq C \gamma\left(\frac{3}{2}-\theta\right)\left(t+T\right)^{-\frac{1}{2}-\theta+\frac{3}{2p}},
}
where $v_i$ is obtained from $f_i$ via the Biot-Savart law.
\end{corollary}
Theorem \ref{theo1} and Corollary \ref{cor1} describe the first order asymptotic profiles of the solutions of the second grade fluids equations. In particular, they show that these solutions behave asymptotically like the self-similar solutions to the heat equation given by
\cent{
$
\displaystyle \left(t,x\right) \longrightarrow \sum^3_{i=1} \frac{b_i}{\left(t+T\right)^2}  f_i\left(\frac{x}{\sqrt{t+T}}\right)
.
$}
In addition, since the same result has been shown in \cite{gallaywayne02bis} for Navier-Stokes equations, it shows that the second grade fluids behave asymptotically like  Newtonian fluids, at least at the first order.
\begin{remark}
We emphasize that the convergence results of Theorem \ref{theo1} and Corollary \ref{cor1} allow to choose the rate of convergence as close as wanted to the optimal one, provided the initial data are small enough in $\mathbb H^2(4)$. In dimension 2, the rate of convergence of the results of Jaffal-Mourtada in \cite{jaffal11} cannot be better than $e^{-\tau/4}$, whereas the optimal one is $e^{-\tau/2}$. In Section \ref{secenergy}, we will see that the method used in the present paper to make estimates on the solutions of (\ref{g2W}) in Sobolev spaces of negative order differs from the one used in \cite{jaffal11}, which is the reason why we are able to obtain a better rate of convergence.
\end{remark}
\begin{remark}
Notice also that the smallness assumption (\ref{inetheo1}) is not optimal. By working harder, we could probably obtain the same theorem with the constant $\gamma\left(\frac{3}{2}-\theta\right)^p$ with $p<2$ in the right hand side of the inequality (\ref{inetheo1}).
\end{remark}
\indent We prove Theorem \ref{theo1} in several steps. First, in Section \ref{secapprox}, we introduce a new system that is close to (\ref{g2W}), but which contains the regularizing term $\varepsilon \Delta^2 W$, with $\varepsilon$ a small positive constant that is devoted to tend to $0$. Due to this regularizing term, we are able, through a semi-group method, to show the existence of local solutions to the regularized system. In a second time, in Section \ref{secenergy} we perform energy estimates on these approximate solutions, and show that these ones are global in time and satisfy the inequality (\ref{inetheo1}). Then, in Section \ref{seclim}, we pass to the limit when $\varepsilon$ tends to $0$ and show that the approximate solutions converge to a global weak solution of (\ref{g2W}) which satisfies the inequality (\ref{inetheo1}). Finally, in order to show that every solution whose initial data satisfy the assumption (\ref{cond1}) converge to his first order asymptotic profile, we show the uniqueness of the weak solutions of (\ref{g2w}) belonging to $C^0\left(\left[0,+\infty\right),\mathbb H^2(4)\right)$.
\paragraph{Biot-Savart law:}
~\\
\vspace{0.2cm}\\
Now, we recall some properties of the Biot-Savart law. Let $w$ be a given divergence free vector field of $\R^3$, the Biot-Savart law gives a divergence free vector field $u$ such that $\curl u=w$. It is given by
\eq{\label{biotsavart}
\displaystyle u(x)=-\frac{1}{4\pi}\int_{\R^3} \frac{\left(x-y\right)\wedge w(y)}{\left|x-y\right|^3} dy.
}
In particular, considering the scaled variables (\ref{scaledvar1}) preserves the Biot-Savart law. Indeed, if $u$ is obtained from $w$ via the Biot-Savart law and $W$ is $w$ expressed into scaled variables, then the divergence free vector field $U$ obtained from $W$ through the Biot-Savart law is $u$ expressed in scaled variables. The next lemma gives some estimates on vector fields obtained by the Biot-Savart law, in various functions spaces.
\lem{\label{biots}
Let $u$ be the velocity field obtained from $w$ via the Biot-Savart law (\ref{biotsavart}).
\begin{enumerate}[(a)]
\item Assume that $1<p<3$, $\frac{3}{2}<q<\infty$ and $\frac{1}{q}=\frac{1}{p}-\frac{1}{3}$. If $w \in L^p(\R^3)^3$, then $u \in L^q(\R^3)^3$, and there exists $C>0$ such that
\eq{\label{biotsa}
\left\|u\right\|_{L^q} \leq C\left\|w\right\|_{L^p}.
}
\item Assume that $1\leq p <3<q\leq \infty$, and define $\eta \in (0,1)$ by the relation $\frac{1}{3}=\frac{\eta}{p}+\frac{\left(1-\eta\right)}{q}$. If $w\in L^p(\R^3)^3\cap L^q(\R^3)^3$, then $u\in L^\infty(\R^3)^3$ and there exists $C>0$ such that
\eq{\label{biotsb}
\left\|u\right\|_{L^\infty}\leq C\left\|w\right\|^\eta_{L^p}\left\|w\right\|^{1-\eta}_{L^q}.
}
\item Assume that $1<p<\infty$. If $w\in L^p(\R^3)^3$, then $\nabla u \in L^p(\R^3)^9$ and there exists $C>0$ such that
\eq{\label{biotsc}
\left\|\nabla u\right\|_{L^p} \leq C\left\|w\right\|_{L^p}.
}
\end{enumerate}
}
This lemma is proved in \cite{gallaywayne02bis} and will be very useful when making estimates on the solutions of (\ref{g2W}).
\section{\label{secapprox}Approximate solutions}
In this section, we introduce a new system that is close to (\ref{g2w}), but contains the regularizing term $\varepsilon \Delta^2 w$, where $\varepsilon$ is a small positive constant. The reason to introduce such a system is to get smooth solutions of the new system, for which we are able to make estimates in $\mathbb H^2(4)$ and obtain the inequality (\ref{inetheo1}). In Section \ref{seclim}, we pass to the limit when $\varepsilon$ goes to $0$ and show that the limit of the solution of the regularized system is a weak solutions of the system (\ref{g2W}) and satisfies also the inequality (\ref{inetheo1}). We introduce the following regularized system, given by
\eq{\label{g2we}
\arr{l}{\displaystyle\partial_t\left(w_\varepsilon-\alpha \Delta w_\varepsilon\right)+\varepsilon \Delta^2 w_\varepsilon -\Delta w_\varepsilon +\curl\left(\left(w_\varepsilon-\alpha \Delta w_\varepsilon\right)\wedge u_\varepsilon\right)=0,\\
\displaystyle\dv u_\varepsilon=\dv w_\varepsilon=0,\\
\displaystyle w_{\varepsilon \left|t=0\right.}=w_0.}
}
The next theorem shows that, for every $w_0 \in \mathbb H^2(4)$, there exists a unique local solution to (\ref{g2we}) belonging to $\mathbb H^2(4)$, which is smooth enough to perform the estimates of Section \ref{secenergy}.
\theo{\label{theoapprox}
Let $\varepsilon >0$ and $w_0 \in \mathbb H^2(4)$. There exists $t_\varepsilon>0$ and a unique solution $w_\varepsilon$ to the system (\ref{g2we}) defined on the time interval $\left[0,t_\varepsilon\right)$ such that
\cent{$
w_\varepsilon \in C^1\left(\left(0,t_\varepsilon\right), \mathbb H^1(4)\right)\cap C^0\left(\left[0,t_\varepsilon\right), \mathbb H^2(4)\right)\cap C^0\left(\left(0,t_\varepsilon\right), \mathbb H^3(4)\right).
$}
}
\textbf{Proof: }To get this result, one defines $w_{\varepsilon,\mu}(t,x) = w_\varepsilon \left(t,\frac{x}{\mu}\right)$, where $\mu>0$. This change of variables enables us to show the existence of solutions to the system (\ref{g2we}) without restrictions on the size of the parameter $\alpha$. We define $u_{\varepsilon,\mu}$ obtained from $w_{\varepsilon,\mu}$ by the Biot-Savart law (\ref{biotsavart}). It is easy to check that $u_{\varepsilon,\mu}(t,x)=\mu u_\varepsilon (t,\frac{x}{\mu})$.   In order to show the existence of a unique solution to (\ref{g2we}), we will prove that there exists a unique solution to the system
\eq{\label{g2ve}
\arr{l}{\partial_t \left(w_{\varepsilon,\mu} - \alpha \mu^2 \Delta w_{\varepsilon,\mu}\right)-\varepsilon \mu^4 \Delta^2 w_{\varepsilon,\mu}- \mu^2\Delta w_{\varepsilon,\mu} + \curl\left(\left(w_{\varepsilon,\mu} - \alpha \mu^2 \Delta w_{\varepsilon,\mu} \right)\wedge u_{\varepsilon,\mu}\right)=0,\\
\dv w_{\varepsilon,\mu} = \dv u_{\varepsilon,\mu}=0,\\
w_{\varepsilon,\mu \left|t=0\right.}=w_0(\frac{x}{\mu}) \in \mathbb H^2(4).

}
}
We define now $\displaystyle z_\varepsilon(t,x)= q(x) w_{\varepsilon,\mu}(t,x)$, where $q(x)=\left(1+\left|x\right|^4\right)$. In particular, if $w_{\varepsilon,\mu} \in \mathbb L^2(4)$, then $z_{\varepsilon} \in \mathcal H$, where
\cent{
$
\mathcal H = \left\{z \in L^2(\R^3)^3 : \dv\left( q^{-1} z\right)=0 \right\}.
$}
For later use, we define, for $s \geq 0$,
\cent{
$
\mathcal H^s = \mathcal H \cap H^s(\R^3)^3, \quad \hbox{and} \quad \mathcal H^{-s} = \left(\mathcal H^s\right)^{'},
$}
where $\left(\mathcal H^s\right)^{'}$ denotes the dual space of $\mathcal H^s$.
\vspace{0.5cm}\\
We equip $\mathcal H^s$ with the classical $H^s$ Sobolev norm, which makes $\mathcal H^s$ complete. From the system (\ref{g2ve}), we deduce the following one, that we solve in $z_\varepsilon$,
\eq{\label{g2ze}
\arr{l}{
\displaystyle \partial_\tau\left(z_\varepsilon-\alpha \mu^2\Delta z_\varepsilon - \alpha \mu^2 q\Delta q^{-1} z_\varepsilon-2\alpha\mu^2 q\nabla q^{-1}.\nabla z_\varepsilon \right)+\varepsilon \mu^4 \Delta^2 z_\varepsilon = F\left(x,z_\varepsilon\right),\\
\dv \left( q^{-1} z_\varepsilon \right)=0,\\
\displaystyle z_{\varepsilon \left|t=0\right.}(x)= z_0 (x) \in \mathcal H^2,
}
}
where
\cent{
$
\arr{l}{
\displaystyle F(x,z_\varepsilon) = -\varepsilon \mu^4 q \Delta ^2 \left(q^{-1} z_\varepsilon\right) +\mu^2 q\Delta \left(q^{-1} z_\varepsilon\right)\\
\hspace{4cm}\displaystyle+q\curl\left(\left(q^{-1} z_\varepsilon -\mu^2 \alpha_1 \Delta\left(q^{-1} z_\varepsilon\right)\right)\wedge u_{\varepsilon,\mu}\right).
}
$
}
The system (\ref{g2ze}) is actually autonomous. Indeed, one can recover $u_{\varepsilon,\mu}$ by the Biot-Savart law (\ref{biotsavart}) applied to $ q^{-1} z_\varepsilon$. To show the existence of solutions to (\ref{g2we}) in $\mathbb H^1(4)$, it suffices to show the existence of solutions to (\ref{g2ze}) in $\mathcal H^1$, for data belonging to $\mathcal H^2$. 
\vspace{0.5cm}\\
We set two linear differential operators $B : D(B)=\mathcal H^1 \rightarrow \mathcal H^{-1}$ and $D : D(D)=\mathcal H \rightarrow \mathcal H^{-1}$, given by
\cent{
$
\arr{l}{B= \alpha \mu^2 q\Delta q^{-1} +\alpha \mu^2 \Delta,\\
D = \alpha \mu^2 q\nabla q^{-1} .\nabla.
}
$
}
Via Lax-Milgram theorem, we show now that if $\mu$ is sufficiently small with respect to $\alpha$, the operator $\left(I-B-D\right)$ is invertible. In order to do that, we define the bilinear form on $\mathcal H^1\times \mathcal H^1$, given by
\cent{
$a(u,v)=\left(u,v\right)_{L^2}+\alpha \mu^2 \left(\nabla u, \nabla v\right)_{L^2} -\alpha \mu^2\left(q\Delta q^{-1}u,v\right)_{L^2}-2\alpha\mu^2\left(q\nabla q^{-1}.\nabla u,v\right)_{L^2}.$
}
Since $q\Delta q^{-1}$ and $q\nabla q^{-1}$ are bounded on $\R^3$, the bilinear form $a$ is continuous on $\mathcal H^1$. We now show, taking $\mu$ small enough, that $a$ is also coercive on $\mathcal H^1$. Indeed, integrating by parts and using Hölder and Young inequalities, we have
\aligne{
a(u,u)&\geq\left(1-\alpha\mu^2 \sup\limits_{x\in \R^3} \left(q\Delta q^{-1}\right)+\alpha \mu^2 \inf\limits_{x\in \R^3} \left(\dv\left(q\nabla q^{-1}\right)\right)\right)\left\|u\right\|^2_{L^2}+\alpha \mu^2 \left\|\nabla u\right\|^2_{L^2}.
}
Thus, if we take $\mu$ sufficiently small, we get
\cent{
$
a(u,u) \geq C(\alpha,\mu)\left\|u\right\|^2_{H^1},
$
}
where $C(\alpha,\mu)$ is a positive constant depending on $\alpha$ and $\mu$.
\vspace{0.5cm}\\
The classical Lax-Milgram theorem enables us to define $\left(I-B-D\right)^{-1}$ from $\mathcal H^{-1}$ to $\mathcal H^1$. We define the linear differential operator $A : D(A)=\mathcal H^3 \rightarrow \mathcal H^1$ given by
\cent{$A= \varepsilon\mu^4 \left(I-B-D\right)^{-1} \Delta^2$.}
We can rewrite the system (\ref{g2ze}) as follows:
\eq{\label{g2zebis}
\arr{l}{
\displaystyle \partial_\tau z_\varepsilon+A z_\varepsilon = \left(I-B-D\right)^{-1}F\left(x,z_\varepsilon\right),\\
\displaystyle z_{\varepsilon \left|t=0\right.}= z_0.
}}
In order to show the existence of solutions to such a system, we use, like in \cite{jaffal11}, a semi-group method. First, we show that $-A$ generates an analytic semi-group on $\mathcal H^1$ which is equivalent as $A$ is sectorial on $\mathcal H^1$.
We decompose $A$ as follows:
\aligne{
A &= \varepsilon \mu^4 \left(Id -B -D\right)^{-1}\Delta^2\\
&=\varepsilon \mu^4 \left(Id -B\right)^{-1}\Delta^2+\varepsilon \mu^4 \left(Id -B -D\right)^{-1} D \left(Id -B  \right)^{-1} \Delta^2\\
&= J+R,
}
where 
\cent{$
\arr{l}{
\displaystyle J=Id+\varepsilon \mu^4 \left(Id -B \right)^{-1}\Delta^2,\\
\displaystyle R=-Id+\varepsilon \mu^4 \left(Id -B -D\right)^{-1} D \left(Id -B  \right)^{-1} \Delta^2.
}
$}
We first show that $J$ is sectorial. We will see later that $R$ satisfies properties that enable to conclude that $A$ is sectorial if $J$ is sectorial. Taking $\mu$ sufficiently small compared to $\alpha$, it is easy, arguing like we did to invert $\left(I-B-D\right)$, to show that $\left(I-B\right)^{-1}$ is well defined from $\mathcal H^{-1}$ to $\mathcal H^1$. Consequently, the operator $J$ is well defined from $\mathcal H^3$ to $\mathcal H^1$. We define now the bilinear form $j$ on $\mathcal H^2\times \mathcal H^2$ associated to $J$. To this end, we introduce a $H^1-$scalar product which is adapted to $J$. We define
\cent{
$
\left\langle u,v \right\rangle_{H^1} = \left(\left(1-\alpha \mu^2 q\Delta q^{-1} \right) u, v\right)_{L^2} + \alpha \mu^2 \left( \nabla u, \nabla v\right)_{L^2}.
$
}
If $\mu$ is sufficiently small, $\left\langle ., . \right\rangle_{H^1}$ is a scalar product on $\mathcal H^1$. In particular, if $u \in \mathcal H^2$ and $v\in \mathcal H^1$, one has
\cent{
$
\left\langle u,v \right\rangle_{H^1} = \left(\left(I-B\right) u,v\right)_{L^2}.
$
}
Via this product, we define
\cent{
$
j(u,v) = \left\langle u,v \right\rangle_{H^1} + \varepsilon \mu^4 \left(\Delta u, \Delta v\right).
$
}
In particular, if $u \in \mathcal H^3$ and $v\in \mathcal H^1$, one has
\cent{
$
j(u,v) = \left\langle Ju,v \right\rangle_{H^1}.
$
}
The bilinear form $j$ is obviously continuous on $\mathcal H^2 \times \mathcal H^2$. Furthermore, if $\mu$ is small enough, it is also coercive on $\mathcal H^2$. Indeed,
\aligne{
j(u,u)& \geq C(\alpha,\mu)\left\|u\right\|^2_{H^1} + \varepsilon \mu^4\left\|\Delta u\right\|^2_{L^2}\\
& \geq C(\alpha, \mu, \varepsilon)\left\|u\right\|^2_{H^2}.
}
Thus $j$ is continuous and coercive on $\mathcal H^2$ and consequently $J$ is sectorial on $\mathcal H^1$, that is equivalent to say that $-J$ generates an analytic semi-group on $\mathcal H^1$. Furthermore, we can check that $R$ is continuous from $\mathcal H^2$ to $\mathcal H^1$, and we have
\cent{
$
\left\|R u\right\|_{H^1} \leq C(\alpha,\mu,\varepsilon )\left\|u\right\|_{H^2}.
$
}
Using the coerciveness of $j$, we get, for all $u \in \mathcal H^3$,
\eq{
\arr{ll}{\left\|R u\right\|^2_{H^1} &\leq C(\alpha,\mu,\varepsilon )j(u,u)\\
&\leq C(\alpha,\mu,\varepsilon ) \left\langle Ju,u \right\rangle_{H^1}\\
&\leq C \left\|J u \right\|_{H^1} \left\|u \right\|_{H^1}.
}
}
Applying the Young inequality, we obtain, for all $\delta>0$
\cent{
$
\left\|R u\right\|^2_{H^1} \leq \delta \left\|J u \right\|^2_{H^1} + C\left\|u \right\|^2_{H^1},
$ for all $ u \in \mathcal H^3$.
}
From a classical result that we can find in the book of D. Henry \cite{henry81}, it implies that $J+R$ is sectorial on $\mathcal H^1$.
\vspace{0.5cm}\\
To achieve this proof, we check that $A^{-1} F(x,v)$ is locally Lipschitz in $v\in \mathcal H^1$ on the bounded sets of $\mathcal H^2$. According to \cite[section 6.3]{pazy83} and \cite[chapter 3]{henry81}, we finally get Theorem \ref{theoapprox}.
\begin{flushright}
$\square$
\end{flushright}
\section{\label{secenergy} Energy estimates}
In this section, we perform several energy estimates on the solution of the system (\ref{g2we}) given by Theorem \ref{theoapprox}. We consider a fixed positive constant $\theta$ such that $0<\theta < \frac{3}{2}$, which is the rate of convergence of Theorem \ref{theo1}. Let $T$ be a positive constant which will be made more precise later and that we assume, without loss of generality, to be such that $T\geq 1$. We consider $W_\varepsilon$ the divergence free vector field obtained from $w_\varepsilon$ via the change of variables (\ref{scaledvar1}). According to Theorem \ref{theoapprox}, there exists a maximal time $\tau_\varepsilon$ such that $W_\varepsilon$ belongs to $C^1\left(\left(\tau_0,\tau_\varepsilon\right),\mathbb H^1(4)\right)\cap C^0\left(\left(\tau_0,\tau_\varepsilon\right),\mathbb H^3(4)\right) $, where $\tau_0= \log (T)$. A short computation shows that $W_\varepsilon$ is the solution of the system
\eq{\label{g2We}
\arr{l}{\partial_\tau\left(W_\varepsilon-\alpha e^{-\tau}\Delta W_\varepsilon\right)+\varepsilon e^{-\tau}\Delta^2 W_\varepsilon - \mathcal L (W_\varepsilon) +\curl\left(\left(W_\varepsilon-\alpha e^{-\tau} \Delta W_\varepsilon\right)\wedge U_\varepsilon\right)\\
\hspace{8cm}+\alpha e^{-\tau}\Delta W_\varepsilon+\alpha e^{-\tau}\frac{X}{2}.\nabla \Delta W_\varepsilon=0,\\
\dv U_\varepsilon=\dv W_\varepsilon=0,\\
W_{\varepsilon \left|\tau=\tau_0\right.}=W_0,}
}
where we recall that
\cent{
$
\mathcal L (W_\varepsilon) = W_\varepsilon+\Delta W_\varepsilon +\frac{X}{2}.\nabla W_\varepsilon.
$}
In this section, we obtain several energy estimates in various functions spaces. More precisely, assuming that $T$ is large enough and $W_0$ is small enough in $\mathbb H^2(4)$, we show that the solution of (\ref{g2We}) stays bounded in time in those energy spaces and is consequently global in time. In addition, we obtain the inequality (\ref{inetheo1}) for $W_\varepsilon$. The method to reach this aim is based on the construction of an energy functional $E$ such that
\cent{
$
\displaystyle E(\tau) \thicksim \left\|W_\varepsilon (\tau)-e^{-\tau}\sum^3_{i=1} b_i f_i\right\|^2_{H^2(4)}
$, for all $\tau \geq \log (T)$,
}
and $E$ satisfies, for all $\tau \geq \log (T)$,
\eq{
\label{ineqenergy}
\partial_\tau E(\tau) + 2\theta E(\tau) \leq C e^{-3\tau},
}
where $\displaystyle b_i = \int_{\R^3} p_i (X) .W_0(X) dX$ and $\left\{f_1,f_2,f_3\right\}$ is the basis of the eigenspace of $\mathcal L$ associated to the eigenvalue $-1$, given by (\ref{fi}). Through the Gronwall Lemma, the inequality (\ref{ineqenergy}) allows to get the inequality (\ref{inetheo1}) for $W_\varepsilon$ and to conclude that $W_\varepsilon$ is global in time.
\vspace{0.5cm}\\
We define $\displaystyle  \Omega_\infty = \sum^3_{i=1} b_i f_i$. The decomposition (\ref{decomp2}) becomes
\eq{\label{decomp3} 
W_\varepsilon (\tau)= e^{-\tau}\Omega_\infty+R_\varepsilon (\tau).
}
A short computation shows that $R_\varepsilon$ satisfies the equality
\eq{\label{g2Re}
\arr{l}{
\partial_\tau\left(R_\varepsilon -\alpha e^{-\tau} \Delta R_\varepsilon \right) +\varepsilon e^{-\tau}\Delta^2 R_\varepsilon -\mathcal L(R_\varepsilon) +\curl\left(\left(W_\varepsilon-\alpha e^{-\tau} \Delta W_\varepsilon\right)\times U_\varepsilon\right) \\
\\
\hspace{4cm}+ \alpha e^{-\tau} \Delta R_\varepsilon+\alpha e^{-\tau}\frac{X}{2}.\nabla \Delta R_\varepsilon +3\alpha e^{-2\tau}\Delta \Omega_\infty+\varepsilon e^{-2\tau} \Delta^2 \Omega_\infty=0.
}
}
In this section, we assume that $W_0$ satisfies the condition (\ref{cond1}) of Theorem \ref{theo1} for some positive constant $\gamma$. We also set $M$ to be a positive constant such that $M\geq 2$ which will be made more precise later. We define $\tau^*_\varepsilon$ the largest positive time such that, for all $\tau \in \left[\tau_0,\tau^*_\varepsilon\right)$,
\eq{\label{cond2W}
\arr{l}{\displaystyle \left\|W_\varepsilon (\tau)\right\|^2_{L^2(4)}+\left\|\nabla W_\varepsilon(\tau)\right\|^2_{L^2}+\alpha e^{-\tau}\left\| \Delta W_\varepsilon(\tau)\right\|^2_{L^2}\\
\displaystyle \hspace{5cm}+\alpha^2 e^{-2\tau}\left\|\left|X\right|^4 \Delta W_\varepsilon(\tau)\right\|^2_{L^2} \leq M \gamma \left(\frac{3}{2}-\theta\right)^2.
}
}
Since $R_\varepsilon$ belongs to $C^0\left(\left[\tau_0,\tau^*_\varepsilon\right),\mathbb H^2(4)\right)$, the time $\tau^*_\varepsilon$ is well defined. The next lemma gives two inequalities on $\nabla W_\varepsilon$ and $R_\varepsilon$.
\lem{\label{conditions}
Let $W_\varepsilon \in C^0\left(\left[\tau_0,\tau^*_\varepsilon\right),\mathbb H^2(4)\right)$ satisfying the condition (\ref{cond2W}) and $R_\varepsilon=W_\varepsilon - e^{-\tau}\Omega_\infty$. There exists a positive constant $C$ such that, for all $\tau \in \left[\tau_0,\tau^*_\varepsilon\right)$,
\eq{\label{cond2}
\arr{l}{\left|b\right|^2+\left\|R_\varepsilon (\tau)\right\|^2_{L^2(4)}+\left\|\nabla R_\varepsilon(\tau)\right\|^2_{L^2}+\alpha e^{-\tau}\left\| \Delta R_\varepsilon(\tau)\right\|^2_{L^2}\\
 \hspace{5cm}+\alpha^2 e^{-2\tau}\left\|\left|X\right|^4 \Delta R_\varepsilon(\tau)\right\|^2_{L^2}\leq C M \gamma\left(\frac{3}{2}-\theta\right)^2.}
}
}
\textbf{Proof: }To prove the inequality (\ref{cond2}), we notice that , for all $i\in \left\{1,2,3\right\}$,
\aligne{
\left|b_i\right| &\leq \int_{\R^2} \left|X\right| \left|W_0\right| dX\\
&\leq \left(\int_{\R^2} \frac{1}{\left(1+\left|X\right|^2\right)^3}dX\right)^{1/2}\left(\int_{\R^2} \left(1+\left|X\right|^2\right)^3 \left|X\right|^2 \left|W_0\right|^2 DX\right)^{1/2}\\
&\leq C \left\|W_0\right\|_{L^2(4)}.
}
Thus, recalling that $\displaystyle R_\varepsilon =W_\varepsilon - e^{-\tau}\sum^3_{i=1} b_i f_i$ and taking into account (\ref{cond1}), we obtain (\ref{cond2}).
\begin{flushright}
$\square$
\end{flushright}
For the sake of simplicity, we assume in this section that $\gamma \leq 1$ and $\left(\frac{3}{2}-\theta\right) \leq 1$.
\subsection{Estimates in $H^{-\left(\theta+2\right)}(\R^3)$}
In this section, we perform an estimate of $R_\varepsilon$ in the space $H^{-\left(\theta+2\right)}(\R^3)$ on the time interval $\left[\tau_0,\tau^*_\varepsilon\right)$. This is motivated by the fact that, in the $H^1-$estimate that we establish below, the term $\left\|R_\varepsilon \right\|^2_{L^2}$ takes place in the right hand side of the inequality (\ref{ineqenergy}). To absorb this term, we look for an estimate in the homogeneous Sobolev space $\dot{H}^{-\left(\theta+2\right)}(\R^3)$. Combined with the other energy estimates, it gives an estimate in the classical Sobolev space $H^{-\left(\theta+2\right)}(\R^3)$. Notice that the constant $\theta + 2$ is chosen in order to obtain the term $2 \theta E$ in the inequality (\ref{ineqenergy}). In \cite{jaffal11}, the choice of the Sobolev space of negative order do not depend on $\theta$, that is why the rate of convergence obtained in \cite{jaffal11} cannot be taken as close as wanted to the optimal one. In order to perform this energy estimate, we define, for $s\in \R$, the operator
\cent{
$\displaystyle \left(-\Delta\right)^{-s} u = \bar{\mathcal F} \left(\frac{1}{\left|\xi\right|^{4s}} \widehat{u}\right)$,
}
where $\widehat{u}$ is the Fourier transform of $u$, given by
\cent{
$\displaystyle \widehat{u}(\xi)=\int_{\R^3} e^{-i x.\xi} u(x) dx,
$
}
and $\bar{\mathcal F}$ is the inverse Fourier transform.
\vspace{0.5cm}\\
In this section, given $0\leq \theta<\frac{3}{2}$, we apply the linear operator $\left(-\Delta\right)^{-\left(\frac{\theta}{2}+1\right)}$ to (\ref{g2Re}) and then make the $L^2-$inner product of it with $\left(-\Delta\right)^{-\left(\frac{\theta}{2}+1\right)} R_\varepsilon$. We are allowed to consider $\left(-\Delta\right)^{-\left(\frac{\theta}{2}+1\right)} R_\varepsilon$ by the lemma
\lem{\label{weight-}
Let $u \in  L^2(4)$ such that $\displaystyle \int_{\R^3} u(x) dx =0$.
\begin{enumerate}
\item  If $\displaystyle \int_{\R^3} x_i u(x) dx =0$ for every $i\in \left\{1,2,3\right\}$, then, for all $0\leq s<\frac{7}{4}$, $\left(-\Delta \right)^{-s} u \in L^2(\R^3)$ and there exists a positive constant $C$ such that
\eq{\label{eqweight-1}
\left\|\left(-\Delta \right)^{-s} u \right\|_{L^2}\leq \frac{C}{\sqrt{7-4s}} \left\| u \right\|_{L^2(4)}.
}
\item For all $0\leq s<\frac{7}{4}$, $\left(-\Delta \right)^{-s} \nabla u \in L^2(\R^3)^3$ and there exists a positive constant $C$ such that
\eq{\label{eqweight-2}
\left\|\left(-\Delta \right)^{-s} \nabla u \right\|_{L^2}\leq \frac{C}{\sqrt{7-4s}} \left\| u \right\|_{L^2(3)}.
}
\end{enumerate}
}
\textbf{Proof: }Using Fourier variables, we get
\aligne{
\left\|\left(-\Delta\right)^{-s} u\right\|^2_{L^2} &=\frac{1}{\left(2\pi\right)^3}\int_{\R^2} \frac{1}{\left|\xi\right|^{4s}}\left|\widehat{u}(\xi)\right|^2 d\xi\\
&\leq \frac{1}{\left(2\pi\right)^3}\int_{\left|\xi\right|\leq 1} \frac{1}{\left|\xi\right|^{4s}}\left|\widehat{u}(\xi)\right|^2 d\xi + \left\|u\right\|^2_{L^2}.
}
We note $\displaystyle I=\frac{1}{\left(2\pi\right)^3}\int_{\left|\xi\right|\leq 1} \frac{1}{\left|\xi\right|^{4s}}\left|\widehat{u}(\xi)\right|^2 d\xi$. Using the fact that $\widehat{u}(0)=\displaystyle \int_{\R^3} u(x) dx =0$ and the Cauchy-Schwartz inequality on the interval $(0,1)$, we have
\aligne{
I&= \frac{1}{\left(2\pi\right)^3}\int_{\left|\xi\right|\leq 1} \frac{1}{\left|\xi\right|^{4s}}\left|\int^1_0 \xi.\nabla \widehat{u}(\sigma\xi) d\sigma\right|^2 d\xi\\
&\leq C\int_{\left|\xi\right|\leq 1} \frac{1}{\left|\xi\right|^{4s-2}}\int^1_0 \left|\nabla \widehat{u}(\sigma\xi) \right|^2 d\sigma d\xi.
}
Then, due to the fact that $\displaystyle\partial_j \widehat{u}(0)= i\int_{\R^2} x_j u(x) dx =0$, we get
\aligne{
I&\leq C \int_{\left|\xi \right|\leq 1} \frac{1}{\left|\xi\right|^{4s-2}}\int^1_0\left(\sum^3_{i,j=1}\left|\int^1_0 \xi_j \partial_i \partial_j \widehat{u}(r \sigma \xi)dr \right|^2\right)  d\sigma d\xi\\
&\leq C\int_{\left|\xi \right|\leq 1} \frac{1}{\left|\xi\right|^{4s-4}}\int^1_0 \int^1_0 \left|\nabla^2 \widehat{u}(r\sigma \xi)\right|^2 dr  d\sigma d\xi.
}
Finally, the continuous injection of $H^2(\R^3)$ into $L^\infty(\R^3)$ yields
\aligne{
I&\leq \frac{C}{7-4s}\left\|\nabla^2 \widehat{u}\right\|^2_{L^\infty}\\
&\leq \frac{C}{7-4s}\left\|\nabla^2 \widehat{u}\right\|^2_{H^2}\\
&\leq \frac{C}{7-4s}\left\|u\right\|^2_{L^2(4)},
}
and thus the inequality (\ref{eqweight-1}) is shown.
\vspace{0.5cm}\\
To get (\ref{eqweight-2}), using Fourier variables, we have
\aligne{
\left\|\left(-\Delta \right)^{-s}\nabla u\right\|^2_{L^2} &= \frac{1}{\left(2\pi\right)^3}\int_{\left|\xi\right|\leq 1} \frac{1}{\left|\xi\right|^{4s-2}} \left|\widehat{u}(\xi)\right|^2 d\xi + \left\|u\right\|^2_{L^2}\\
&= \frac{1}{\left(2\pi\right)^3}\int_{\left|\xi\right|\leq 1} \frac{1}{\left|\xi\right|^{4s-2}} \left|\int^1_0 \xi.\nabla\widehat{u}(s\xi)ds\right|^2 d\xi + \left\|u\right\|^2_{L^2}\\
&\leq \frac{1}{\left(2\pi\right)^3}\int_{\left|\xi\right|\leq 1} \frac{1}{\left|\xi\right|^{4s-4}} \left|\int^1_0 \left|\nabla\widehat{u}(s\xi)\right|ds\right|^2 d\xi + \left\|u\right\|^2_{L^2}.
}
Using now Hölder inequalities, the fact that $4s-4 <3$ and the continuous injection of $H^2(\R^3)$ into $L^\infty(\R^3)$, we have
\aligne{
\left\|\left(-\Delta \right)^{-s}\nabla u\right\|^2_{L^2} &\leq C\int^1_0\int_{\left|\xi\right|\leq 1} \frac{1}{\left|\xi\right|^{4s-4}} \left|\nabla\widehat{u}(s\xi)\right|^2 d\xi ds+ \left\|u\right\|^2_{L^2}\\
&\leq C\left(\int_{\left|\xi\right|\leq 1} \frac{1}{\left|\xi\right|^{4s-4}}d\xi\right) \left\|\nabla \widehat{u}\right\|^2_{L^\infty}+ \left\|u\right\|^2_{L^2}\\
&\leq \frac{C}{7-4s} \left\| u\right\|^2_{L^2(3)}+ \left\|u\right\|^2_{L^2}.
}
\begin{flushright}
$\square$
\end{flushright}
In order to apply the lemma \ref{weight-} to the non linear terms of the equation (\ref{g2Re}), we state the following lemma.
\lem{\label{gnouf2}
Let $w \in \mathbb H^2(4)$ and $u$ obtained from $w$ via the Biot-Savart law (\ref{biotsavart}). For all $C\in \R$, we have
\eq{
\int_{\R^3} \left(w(x)-C\Delta w(x)\right)\wedge u(x) dx=0.
}
}
\textbf{Proof: }In order to show this equality, we just have to look at the equality (\ref{gnouf}). An integration by parts gives directly (\ref{gnouf2}).
\begin{flushright}
$\square$
\end{flushright}
\lem{\label{H-L}
Let $w$ belongs to $H^2(4)$ and $s$ such that $0\leq s <\frac{7}{4}$, then $u$ satisfies the equalities
\begin{enumerate}
\item $\left(\left(-\Delta\right)^{-s}\mathcal{L}(w),\left(-\Delta\right)^{-s} w\right)_{L^2} = -\left\|\left(-\Delta \right)^{\frac{1}{2}-s} w\right\|^2_{L^2}-\left(s-\frac{1}{4}\right) \left\|\left(-\Delta \right)^{-s} w\right\|^2_{L^2}$.
\item $\left(\left(-\Delta\right)^{-s}\left(\frac{x}{2}.\nabla \Delta w\right),\left(-\Delta\right)^{-s} w\right)_{L^2} = \left(s+\frac{5}{4}\right)\left\|\left(-\Delta \right)^{\frac{1}{2}-s} w\right\|^2_{L^2}$.
\end{enumerate}
}
This lemma is easily obtained with a few integrations by parts, when passing into Fourier variables.
\vspace{0.5cm}\\
\indent In this section, to simplify the notations, we note $R$ instead of $R_\varepsilon$, $W$ instead of $W_\varepsilon$ and $U$ instead of $U_\varepsilon$. We also note $V_\infty$, the divergence free vector field obtained from $\Omega_\infty$ via the Biot-Savart law and $K$ the divergence free vector field obtained from $R$ via the Biot-Savart law. We assume also, without loss of generality, that $T$ is sufficiently large so that $\alpha e^{-\tau_0}\leq 1$, where we recall that $\tau_0 = \log(T)$. We define the energy functional
\cent{
$
\displaystyle E_0(\tau)=\frac{1}{2}\left(\left\|\left(-\Delta\right)^{-\left(\frac{\theta}{2}+1\right)}R\right\|^2_{L^2}+\alpha e^{-\tau}\left\|\left(-\Delta\right)^{-\left(\frac{\theta+1}{2}\right)}R\right\|^2_{L^2}\right).
$
}
The next lemma gives a $H^{-\left(\theta+2\right)}$ which is necessary to obtain a good rate of convergence in Theorem \ref{theo1}.
\lem{\label{lemE0}
Let $W \in C^1\left(\left(\tau_0,\tau_\varepsilon\right),\mathbb H^1(4)\right)\cap  C^0\left(\left(\tau_0,\tau_\varepsilon\right),\mathbb H^3(4)\right)$ be the solution of (\ref{g2We}). There exist two positive constant $\gamma_0$ and $T_0$ such that, if $T \geq T_0$ and $W_\varepsilon$ satisfies the condition (\ref{cond2W}) for some $\gamma$ such that $0<\gamma \leq \gamma_0$, then there exists a positive constant $C$ such that, for all $\tau \in \left[\tau_0, \tau^*_\varepsilon\right)$,
\eq{\label{eqE0}
\arr{l}{\displaystyle \partial_\tau E_0 +2\theta E_0+\frac{1}{2}\left\|\left(-\Delta \right)^{-\left(\frac{\theta+1}{2}\right)} R\right\|^2_{L^2} \leq \\
\hspace{1cm}\displaystyle C M \gamma \left(\left\|\left|X\right|^4 R\right\|^2_{L^2}+\left\|\nabla R\right\|^2_{L^2}+\alpha^2 e^{-2\tau} \left\|\Delta R \right\|^2_{L^2(4)}\right)+C M^2\gamma\left(\frac{3}{2}-\theta\right) e^{-4\tau}.}
}
}
\textbf{Proof: }To prove this lemma, we apply the operator $\left(-\Delta\right)^{-\left(\frac{\theta}{2}+1\right)}$ to (\ref{g2Re}) and make the $L^2-$inner product of it with $\left(-\Delta\right)^{-\left(\frac{\theta}{2}+1\right)} R$. Applying Lemma \ref{H-L} and through some easy computations, one has
\eq{\label{eqE01}\arr{l}{
\displaystyle \frac{1}{2}\partial_\tau\left(\left\|\left(-\Delta \right)^{-\left(\frac{\theta}{2}+1\right)}R\right\|^2_{L^2} +\alpha e^{-\tau}\left\|\left(-\Delta \right)^{-\left(\frac{\theta+1}{2}\right)}R\right\|^2_{L^2}\right)+\varepsilon e^{-\tau} \left\|\left(-\Delta\right)^{-\frac{\theta}{2} } R\right\|^2_{L^2}\\
\\
\displaystyle \hspace{0.1cm}+\left(\frac{\theta}{2}+\frac{3}{4}\right)\left\|\left(-\Delta \right)^{-\left(\frac{\theta}{2}+1\right)} R\right\|^2_{L^2}+\left(1+ \left(\frac{\theta}{2}+\frac{3}{4}\right) \alpha e^{-\tau}\right)\left\|\left(-\Delta \right)^{-\left(\frac{\theta+1}{2}\right)} R\right\|^2_{L^2} = I_1+I_2,
}
}
where
\cent{$
\arr{l}{
\displaystyle I_1 = \left(\left(-\Delta\right)^{-\left(\frac{\theta}{2}+1\right)}\left(\curl\left(\left( W-\alpha e^{-\tau} \Delta W\right)\wedge U\right)\right), \left(-\Delta\right)^{-\left(\frac{\theta}{2}+1\right)} R\right)_{L^2},\\
\\
\displaystyle I_2 =e^{-2\tau}\left(\left(-\Delta\right)^{-\left(\frac{\theta}{2}+1\right)}\left(-\alpha \Delta \Omega_\infty-\varepsilon \Delta^2 \Omega_\infty\right), \left(-\Delta\right)^{-\left(\frac{\theta}{2}+1\right)} R\right)_{L^2}.
}
$}
We start with the estimate of the easiest term, that is $I_2$. Using the Cauchy-Schwartz inequality, we get
\aligne{
I_2 \leq \alpha e^{-2\tau} \left\|\left(-\Delta\right)^{-\frac{\theta}{2}} \Omega_\infty\right\|_{L^2}&\left\|\left(-\Delta\right)^{-\left(\frac{\theta}{2}+1\right)} R\right\|_{L^2} \\
&+ \varepsilon e^{-2\tau} \left\|\left(-\Delta\right)^{1-\frac{\theta}{2}} \Omega_\infty\right\|_{L^2}\left\|\left(-\Delta\right)^{-\left(\frac{\theta}{2}+1\right)} R\right\|_{L^2}.
}
Using the Lemma \ref{weight-}, the Young inequality and taking into account the good regularity of $\Omega_\infty$ and the inequality (\ref{cond2}), one has
\eq{\label{eqE0I2}\arr{ll}{
I_2 & \displaystyle \leq C e^{-2\tau}\left\|\Omega_\infty\right\|_{H^2(4)}\left\|\left(-\Delta\right)^{-\left(\frac{\theta}{2}+1\right)} R\right\|_{L^2}\\
&\displaystyle \leq \mu \left\|\left(-\Delta\right)^{-\left(\frac{\theta}{2}+1\right)} R\right\|^2_{L^2}+\frac{C\left|b\right|^2}{\mu}e^{-4\tau}\\
&\displaystyle \leq \mu \left\|\left(-\Delta\right)^{-\left(\frac{\theta}{2}+1\right)} R\right\|^2_{L^2}+\frac{CM \gamma \left(\frac{3}{2}-\theta\right)^2}{\mu}e^{-4\tau},\\
}
}
where $\mu$ is a positive constant that will be made more precise later.
\vspace{0.5cm}\\
It remains to bound $I_1$. Using the Cauchy-Schwartz inequality and the lemmas \ref{gnouf2} and \ref{weight-}, we obtain
\aligne{
I_1 &\leq C\left\|\left(-\Delta\right)^{-\left(\frac{\theta}{2}+1\right)}\nabla \left(\left( W-\alpha e^{-\tau} \Delta W\right)\wedge U\right)\right\|_{L^2}\left\|\left(-\Delta\right)^{-\left(\frac{\theta}{2}+1\right)} R\right\|_{L^2}\\
&\leq \frac{C}{\left(\frac{3}{2}-\theta\right)^{1/2}}\left\|\left( W-\alpha e^{-\tau} \Delta W\right)U\right\|_{L^2(4)}\left\|\left(-\Delta\right)^{-\left(\frac{\theta}{2}+1\right)} R\right\|_{L^2}\\
&\leq \frac{C}{\left(\frac{3}{2}-\theta\right)^{1/2}}\left\|U\right\|_{L^\infty}\left\| W-\alpha e^{-\tau} \Delta W\right\|_{L^2(4)}\left\|\left(-\Delta\right)^{-\left(\frac{\theta}{2}+1\right)} R\right\|_{L^2}.\\
}
The inequality (\ref{biotsb}) of Lemma \ref{biots} with $p=2$, $q=6$ and $\eta=\frac{1}{2}$ and the continuous injection of $H^1(\R^3)$ into $L^6(\R^3)$ yield
\aligne{
I_1 &\leq \frac{C}{\left(\frac{3}{2}-\theta\right)^{1/2}}\left\|W\right\|^{1/2}_{L^2}\left\|W\right\|^{1/2}_{L^6}\left\| W-\alpha e^{-\tau} \Delta W\right\|_{L^2(4)}\left\|\left(-\Delta\right)^{-\left(\frac{\theta}{2}+1\right)} R\right\|_{L^2}\\
&\leq \frac{C}{\left(\frac{3}{2}-\theta\right)^{1/2}} \left\|W\right\|_{H^1}\left(\left\| W\right\|_{L^2(4)}+\alpha e^{-\tau} \left\|\Delta W\right\|_{L^2(4)}\right)\left\|\left(-\Delta\right)^{-\left(\frac{\theta}{2}+1\right)} R\right\|_{L^2}\\
&\leq \mu \left\|\left(-\Delta\right)^{-\left(\frac{\theta}{2}+1\right)} R\right\|^2_{L^2}\\
&\hspace{2cm}+\frac{C }{\mu\left(\frac{3}{2}-\theta\right)}\left(\left\|W\right\|^2_{L^2}+\left\|\nabla W\right\|^2_{L^2}\right)\left(\left\| W\right\|^2_{L^2(4)}+\alpha^2 e^{-2\tau} \left\|\Delta W\right\|^2_{L^2(4)}\right).
}
Due to the decomposition (\ref{decomp3}), one has
\aligne{
I_1 &\leq \mu \left\|\left(-\Delta\right)^{-\left(\frac{\theta}{2}+1\right)} R\right\|^2_{L^2}\\
&\hspace{0.5cm}+\frac{C }{\mu\left(\frac{3}{2}-\theta\right)}\left(\left\|R\right\|^2_{L^2}+\left\|\nabla R\right\|^2_{L^2}\right)\left(\left\| W\right\|^2_{L^2(4)}+\alpha^2 e^{-2\tau} \left\|\Delta W\right\|^2_{L^2(4)}\right)\\
&\hspace{1cm}+\frac{C e^{-2\tau}}{\mu\left(\frac{3}{2}-\theta\right)}\left(\left\|\Omega_\infty \right\|^2_{L^2}+\left\|\nabla \Omega_\infty\right\|^2_{L^2}\right)\left(\left\| R \right\|^2_{L^2(4)}+\alpha^2 e^{-2\tau} \left\|\Delta R \right\|^2_{L^2(4)}\right)\\
&\hspace{1.5cm}+\frac{C e^{-4\tau}}{\mu\left(\frac{3}{2}-\theta\right)}\left(\left\|\Omega_\infty \right\|^2_{L^2}+\left\|\nabla \Omega_\infty\right\|^2_{L^2}\right)\left(\left\| \Omega_\infty \right\|^2_{L^2(4)}+\alpha^2 e^{-2\tau} \left\|\Delta \Omega_\infty \right\|^2_{L^2(4)}\right).\\
}
Finally, using the inequalities (\ref{cond2W}) and (\ref{cond2}), we obtain
\eq{\label{eqE0I1}
\arr{l}{\displaystyle I_1 \leq \mu \left\|\left(-\Delta\right)^{-\left(\frac{\theta}{2}+1\right)} R\right\|^2_{L^2}+\frac{C M^2\gamma^2 \left(\frac{3}{2}-\theta\right)^3 e^{-4\tau}}{\mu}\\
\displaystyle \hspace{2cm}+\frac{C M \gamma\left(\frac{3}{2}-\theta\right)}{\mu}\left(\left\|R\right\|^2_{L^2(4)}+\left\|\nabla R\right\|^2_{L^2}+\alpha^2 e^{-2\tau} \left\|\Delta R \right\|^2_{L^2(4)}\right).
}
}
Combining (\ref{eqE01}), (\ref{eqE0I2}) and (\ref{eqE0I1}), it comes
\eq{\arr{l}{
\displaystyle \frac{1}{2}\partial_\tau\left(\left\|\left(-\Delta \right)^{-\left(\frac{\theta}{2}+1\right)}R\right\|^2_{L^2} +\alpha e^{-\tau}\left\|\left(-\Delta \right)^{-\left(\frac{\theta+1}{2}\right)}R\right\|^2_{L^2}\right)+\varepsilon e^{-\tau} \left\|\left(-\Delta\right)^{-\frac{\theta}{2} } R\right\|^2_{L^2}\\
\displaystyle \hspace{6.5cm}+\left(\theta+\frac{1}{2}\left(\frac{3}{2}-\theta-2\mu\right)\right)\left\|\left(-\Delta \right)^{-\left(\frac{\theta}{2}+1\right)} R\right\|^2_{L^2}\\
\displaystyle \hspace{7.5cm}+\left(1+ \left(\frac{\theta}{2}+\frac{3}{4}\right) \alpha e^{-\tau}\right)\left\|\left(-\Delta \right)^{-\left(\frac{\theta+1}{2}\right)} R\right\|^2_{L^2} \\
\\
\displaystyle \hspace{0.5cm}\leq \frac{C M \gamma\left(\frac{3}{2}-\theta\right)}{\mu}\left(\left\|R\right\|^2_{L^2(4)}+\left\|\nabla R\right\|^2_{L^2}+\alpha^2 e^{-2\tau} \left\|\Delta R \right\|^2_{L^2(4)}\right)+\frac{C M^2\gamma\left(\frac{3}{2}-\theta\right)^2 e^{-4\tau}}{\mu}.
}
}
We set $\displaystyle \mu = \frac{\frac{3}{2}-\theta}{2}$, and we obtain
\eq{\label{eqE02}\arr{l}{
\displaystyle \partial_\tau E_0 +2\theta E_0+\left\|\left(-\Delta \right)^{-\left(\frac{\theta+1}{2}\right)} R\right\|^2_{L^2} \leq \\
\hspace{2cm}\displaystyle  C M \gamma\left(\left\|R\right\|^2_{L^2(4)}+\left\|\nabla R\right\|^2_{L^2}+\alpha^2 e^{-2\tau} \left\|\Delta R \right\|^2_{L^2(4)}\right)+C M^2\gamma \left(\frac{3}{2}-\theta\right) e^{-4\tau}.
}
}
Furthermore, using Fourier variables and Hölder inequalities, we see that
\aligne{
\left\|R\right\|^2_{L^2}&=\frac{1}{\left(2\pi\right)^3}\int_{\R^3} \left|\widehat{R}(\xi)\right|^2 d\xi\\
& \leq \frac{1}{\left(2\pi\right)^3}\int_{\R^3} \left|\xi\right|^{\frac{2\left(1+\theta\right)}{2+\theta}}\left|\widehat{R}(\xi)\right|^{\frac{2\left(1+\theta\right)}{2+\theta}} \frac{1}{\left|\xi\right|^{\frac{2\left(1+\theta\right)}{2+\theta}}}\left|\widehat{R}(\xi)\right|^{\frac{2}{2+\theta}} d\xi\\
&\leq \left(\frac{1}{\left(2\pi\right)^3}\int_{\R^3}\frac{1}{\left|\xi \right|^{2\left(\theta+1\right)}}\left|\widehat{R}(\xi)\right|^2 d\xi\right)^{\frac{1+\theta}{2+\theta}}\left(\frac{1}{\left(2\pi\right)^3}\int_{\R^3}\left|\xi \right|^2\left|\widehat{R}(\xi)\right|^2 d\xi\right)^{\frac{1}{2+\theta}}\\
&\leq \left\|\left(-\Delta \right)^{-\left(\frac{\theta+1}{2}\right)} R\right\|^{\frac{2\left(1+\theta\right)}{2+\theta}}_{L^2}\left\|\nabla R\right\|^{\frac{2}{2+\theta}}_{L^2}.
}
Using a convexity inequality, it is easy to see that
\cent{
$\displaystyle\left\|R\right\|^2_{L^2} \leq \frac{1}{\eta^{\frac{2+\theta}{1+\theta}}} \left(\frac{1+\theta}{2+\theta}\right)\left\|\left(-\Delta \right)^{-\left(\frac{\theta+1}{2}\right)} R\right\|^2_{L^2}+ \frac{\eta^{2+\theta}}{2+\theta} \left\|\nabla R\right\|^2_{L^2},$
}
for all $0<\eta\leq 1$.
\vspace{0.5cm}\\
Via a short computation, using the fact that $0<\theta<\frac{3}{2}$ and $0<\eta\leq 1$, we obtain
\eq{\label{interp}
\left\|R\right\|^2_{L^2} \leq \frac{5}{7\eta^2}\left\|\left(-\Delta \right)^{-\left(\frac{\theta+1}{2}\right)} R\right\|^2_{L^2}+ \frac{\eta^2}{2} \left\|\nabla R\right\|^2_{L^2}.
}
Applying (\ref{interp}) with $\eta=1$ and taking $\gamma$ small enough, the inequality (\ref{eqE02}) becomes
\eq{\label{eqE03}
\arr{l}{\displaystyle \partial_\tau E_0 +2\theta E_0+\frac{1}{2}\left\|\left(-\Delta \right)^{-\left(\frac{\theta+1}{2}\right)} R\right\|^2_{L^2} \leq \\
\hspace{1cm}\displaystyle C M \gamma  \left(\left\|\left|X\right|^4 R\right\|^2_{L^2}+\left\|\nabla R\right\|^2_{L^2}+\alpha^2 e^{-2\tau} \left\|\Delta R \right\|^2_{L^2(4)}\right)+C M^2\gamma \left(\frac{3}{2}-\theta\right) e^{-4\tau}.}
}
\begin{flushright}
$\square$
\end{flushright}
\subsection{\label{secH1}Estimates in $H^1(\R^3)$}
This section is devoted to the $H^1-$estimate of the solutions of (\ref{g2Re}) under the condition (\ref{cond2W}). In particular, we see in this section that the previous estimate in $\dot{H}^{-\left(1+\theta\right)}$ enables to absorb the terms involving the $L^2-$norm of $R$. To obtain this $H^1-$estimate, we make the $L^2-$scalar product of (\ref{g2Re}) with $R$. We define the energy functional
\cent{
$
\displaystyle E_1(\tau) = \frac{1}{2}\left(\left\|R\right\|^2_{L^2} + \alpha e^{-\tau} \left\|\nabla R\right\|^2_{L^2}\right).
$
}
The estimate of $R$ in the Sobolev space $H^1(\R^3)$ is given by the next lemma.
\lem{\label{lemE1}
Let $W\in C^1\left(\left(\tau_0, \tau_\varepsilon\right),\mathbb H^1(4)\right)\cap C^0\left(\left(\tau_0, \tau_\varepsilon\right),\mathbb H^3(4)\right)$ be the solution of (\ref{g2We}). There exist two positive constants $\gamma_0$ and $T_0$ such that, if $T \geq T_0$ and $W$ satisfies the condition (\ref{cond2W}) for some $\gamma$ such that $0<\gamma\leq \gamma_0$, then there exists a positive constant $C$ such that, for all $\tau \in \left[\tau_0,\tau^*_\varepsilon\right)$,
\eq{\label{eqE1}
\arr{l}{\displaystyle\partial_\tau E_1 + 3 E_1 +\frac{1}{2}\left\|\nabla R\right\|^2_{L^2}\leq \frac{7}{4}\left\|R\right\|^2_{L^2}+ C M^2 \gamma \left(\frac{3}{2}-\theta\right)^2 e^{-4\tau}\\
\displaystyle \hspace{5cm} +C M\gamma \left(\frac{3}{2}-\theta\right)^2 \left(\left\|R\right\|^2_{L^2}+\alpha^2 e^{-2\tau}\left\|\Delta R\right\|^2_{L^2}\right).
}
}}
\textbf{Proof: }We perform the $L^2-$scalar product of (\ref{g2Re}) with $R$. Performing several integrations by parts, we obtain
\eq{\label{E11}
\frac{1}{2}\partial_\tau \left(\left\|R\right\|^2_{L^2}+\alpha e^{-\tau} \left\|\nabla R\right\|^2_{L^2}\right)+\varepsilon\left\|\Delta R\right\|^2_{L^2}+\left(1-\frac{\alpha}{4}e^{-\tau}\right)\left\|\nabla R\right\|^2_{L^2}-\frac{1}{4}\left\|R\right\|^2_{L^2}=I_1+I_2,
}
where
\cent{$
\arr{l}{
I_1 = \left(\curl\left(\left(W-\alpha e^{-\tau}\Delta W\right)\wedge U\right), R\right)_{L^2},\\
\\
I_2= e^{-2\tau}\left(-\alpha \Delta \Omega_\infty-\varepsilon\Delta^2 \Omega_\infty, R\right)_{L^2}.
}$
}
As usual, because of the good regularity of $\Omega_\infty$, the easiest term to estimate is $I_2$. Integrating by parts, one has
\aligne{
I_2 & = e^{-2\tau}\left(\alpha \nabla \Omega_\infty+\varepsilon\nabla \Delta \Omega_\infty, \nabla R\right)_{L^2}.
}
Using the Hölder and Young inequalities and the inequality (\ref{cond2}), we get
\eq{\label{eqE1I2}
\arr{ll}{
\displaystyle I_2 &\displaystyle\leq e^{-2\tau} \left(\alpha \left\|\nabla \Omega_\infty\right\|_{L^2}+\varepsilon \left\|\nabla \Delta \Omega_\infty\right\|_{L^2}\right)\left\|\nabla R\right\|_{L^2}\\
& \displaystyle\leq C \left|b\right|\left(\alpha+\varepsilon\right)e^{-2\tau}\left\|\nabla R\right\|_{L^2}\\
&\displaystyle\leq \mu \left\|\nabla R\right\|^2_{L^2}+\frac{CM\gamma \left(\frac{3}{2}-\theta\right)^2}{\mu}e^{-4\tau},
}
}
where $\mu$ is a positive constant that will be made more precise later.
\vspace{0.5cm}\\
The last remaining term will be estimated by the same way, using the divergence free property of $U$. Integrating by parts, we obtain
\aligne{
I_1 = \left(\left(W-\alpha e^{-\tau}\Delta W\right)\wedge U,\curl R\right)_{L^2}.
}
We recall that $\curl K = R$ and $\curl V_\infty = \Omega_\infty$ and we decompose $I_1$ as the sum of three terms
\cent{$
I_1=I^1_1+I^2_1+I^3_1,
$}
where
\cent{
$
\arr{l}{
\displaystyle I^1_1 = \left(\left(W-\alpha e^{-\tau}\Delta W\right)\wedge K,\curl R\right)_{L^2},\\
\\
\displaystyle I^2_1 = e^{-\tau} \left( \left(R-\alpha e^{-\tau}\Delta R\right)\wedge V_\infty,\curl R\right)_{L^2},\\
\\
\displaystyle I^3_1 = e^{-2\tau}\left(\left(\Omega_\infty-\alpha e^{-\tau}\Delta \Omega_\infty\right)\wedge V_\infty,\curl R\right)_{L^2}.
}
$}
The Hölder inequalities lead to
\aligne{
I^1_1 & \leq C\left(\left\|K W\right\|_{L^2} +\alpha e^{-\tau} \left\|K\Delta W\right\|_{L^2}\right)\left\|\nabla R\right\|_{L^2}\\
&\leq C\left\|K\right\|_{L^\infty}\left(\left\|W\right\|_{L^2} +\alpha e^{-\tau} \left\|\Delta W\right\|_{L^2}\right)\left\|\nabla R\right\|_{L^2}.
}
Applying the inequality (\ref{biotsb}) with $p=2$, $q=6$ and $\eta=\frac{1}{2}$ and using the continuous injection of $H^1(\R^3)$ into $L^6(\R^3)$, one gets
\aligne{
I^1_1 &\leq C\left\|R\right\|^{1/2}_{L^2}\left\|R\right\|^{1/2}_{L^6}\left(\left\|W\right\|_{L^2} +\alpha e^{-\tau} \left\|\Delta W\right\|_{L^2}\right)\left\|\nabla R\right\|_{L^2}\\
&\leq C\left\|R\right\|^{1/2}_{L^2}\left\|R\right\|^{1/2}_{H^1}\left(\left\|W\right\|_{L^2} +\alpha e^{-\tau} \left\|\Delta W\right\|_{L^2}\right)\left\|\nabla R\right\|_{L^2}.
}
Then, we use the Young inequality and the inequality (\ref{cond2W}). We obtain
\aligne{
I^1_1 &\leq \mu \left\|\nabla R\right\|^2_{L^2} +\frac{C}{\mu}\left(\left\|W\right\|^2_{L^2} +\alpha^2 e^{-2\tau} \left\|\Delta W\right\|^2_{L^2}\right)\left(\left\|R\right\|^2_{L^2}+\left\|\nabla R\right\|^2_{L^2}\right)\\
&\leq \mu \left\|\nabla R\right\|^2_{L^2} +\frac{C M \gamma \left(\frac{3}{2}-\theta\right)^2}{\mu}\left(\left\|R\right\|^2_{L^2}+\left\|\nabla R\right\|^2_{L^2}\right).
}
The Hölder inequalities yield
\aligne{
I^2_1 &\leq C e^{-\tau} \left\|V_\infty\right\|_{L^\infty}\left(\left\|R\right\|_{L^2}+\alpha e^{-\tau}\left\|\Delta R\right\|_{L^2}\right)\left\|\nabla R\right\|_{L^2}.
}
Applying the inequality (\ref{biotsb}) of the lemma \ref{biots} with $p=2$, $q=6$ and $\eta =\frac{1}{2}$, and the inequality (\ref{cond2}), we get
\aligne{
I^2_1 &\leq Ce^{-\tau} \left\|\Omega_\infty\right\|^{1/2}_{L^2}\left\|\Omega_\infty\right\|^{1/2}_{L^6}\left(\left\|R\right\|_{L^2}+\alpha e^{-\tau}\left\|\Delta R\right\|_{L^2}\right)\left\|\nabla R\right\|_{L^2}\\
&\leq C \left|b\right|e^{-\tau}  \left(\left\|R\right\|_{L^2}+\alpha e^{-\tau}\left\|\Delta R\right\|_{L^2}\right)\left\|\nabla R\right\|_{L^2}\\
&\leq \mu \left\|\nabla R\right\|^2_{L^2}+ \frac{C M\gamma \left(\frac{3}{2}-\theta\right)^2}{\mu} \left(\left\|R\right\|^2_{L^2}+\alpha^2 e^{-2\tau}\left\|\Delta R\right\|^2_{L^2}\right).
}
It remains to estimate $I^3_1$. By the same computations, we get
\aligne{
I^3_1 &\leq \mu \left\|\nabla R\right\|^2_{L^2}+\frac{C}{\mu}e^{-4\tau}\left\|V_\infty\right\|^2_{L^\infty}\left(\left\|\Omega_\infty\right\|^2_{L^2}+\alpha^2 e^{-2\tau} \left\|\Delta \Omega_\infty\right\|^2_{L^2}\right)\\
&\leq \mu \left\|\nabla R\right\|^2_{L^2}+\frac{C}{\mu}e^{-4\tau}\left\|\Omega_\infty\right\|_{L^2}\left\|\Omega_\infty\right\|_{L^6}\left(\left\|\Omega_\infty\right\|^2_{L^2}+\alpha^2 e^{-2\tau} \left\|\Delta \Omega_\infty\right\|^2_{L^2}\right)\\
&\leq \mu \left\|\nabla R\right\|^2_{L^2}+\frac{CM^2 \gamma^2 \left(\frac{3}{2}-\theta\right)^4}{\mu}e^{-4\tau}.
}
In particular, we have shown that
\eq{\label{eqE1I1}
\arr{l}{\displaystyle I_1 \leq 3\mu \left\|\nabla R\right\|^2_{L^2} + \frac{C M^2 \gamma^2 \left(\frac{3}{2}-\theta\right)^4}{\mu}e^{-4\tau}\\
\hspace{1cm}\displaystyle + \frac{CM\gamma \left(\frac{3}{2}-\theta\right)^2}{\mu}\left(\left\|R\right\|^2_{L^2}+\left\|\nabla R\right\|^2_{L^2}+\alpha^2 e^{-2\tau}\left\|\Delta R\right\|^2_{L^2}\right).
}
}
Thus, due to the inequalities (\ref{eqE1I2}) and (\ref{eqE1I1}), the inequality (\ref{E11}) becomes
\eq{\label{E12}
\arr{l}{\displaystyle\partial_\tau E_1 + 3 E_1 +\left(1-4\mu-\frac{7\alpha}{4}e^{-\tau}\right)\left\|\nabla R\right\|^2_{L^2}\leq \frac{7}{4}\left\|R\right\|^2_{L^2}+ \frac{CM^2 \gamma \left(\frac{3}{2}-\theta\right)^2}{\mu}e^{-4\tau} \\
\\
\displaystyle\hspace{5cm}+ \frac{C M\gamma \left(\frac{3}{2}-\theta\right)^2}{\mu}\left(\left\|R\right\|^2_{L^2}+\left\|\nabla R\right\|^2_{L^2}+\alpha^2 e^{-2\tau}\left\|\Delta R\right\|^2_{L^2}\right).
}
}
Taking $\gamma_0$ and $\mu$ small enough and $T = e^{\tau_0}$ large enough, we obtain the inequality
\eq{\label{E13}
\arr{l}{\displaystyle\partial_\tau E_1 + 3 E_1 +\frac{1}{2}\left\|\nabla R\right\|^2_{L^2}\leq \frac{7}{4}\left\|R\right\|^2_{L^2}+ C M\gamma \left(\frac{3}{2}-\theta\right)^2 \left(\left\|R\right\|^2_{L^2}+\alpha^2 e^{-2\tau}\left\|\Delta R\right\|^2_{L^2}\right)\\
\displaystyle \hspace{11cm} +C M^2 \gamma \left(\frac{3}{2}-\theta\right)^2 e^{-4\tau},
}
}
that concludes the proof of this lemma.
\begin{flushright}
$\square$
\end{flushright}
In order to achieve the $H^1-$estimate of $R$, we now combine the energy inequalities (\ref{eqE0}) and (\ref{eqE1}). Using the interpolation inequality (\ref{interp}), we get, from the inequality (\ref{eqE1}),
\eq{
\arr{l}{\displaystyle\partial_\tau E_1 + 3 E_1 +\frac{1}{2}\left\|\nabla R\right\|^2_{L^2}\leq \frac{7}{4}\left(\frac{5}{7\eta^2}\left\|\left(-\Delta\right)^{-\left(\frac{1+\theta}{2}\right)}R\right\|^2_{L^2}+\frac{\eta^2}{2}\left\|\nabla R\right\|^2_{L^2}\right)\\
\hspace{1cm}+ C M\gamma \left(\frac{3}{2}-\theta\right)^2\left(\left\|\left(-\Delta\right)^{-\left(\frac{1+\theta}{2}\right)}R\right\|^2_{L^2}+\left\|\nabla R\right\|^2_{L^2}+\alpha^2 e^{-2\tau}\left\|\Delta R\right\|^2_{L^2}\right)\\
\displaystyle \hspace{10cm}+ C M^2 \gamma \left(\frac{3}{2}-\theta\right)^2 e^{-4\tau},
}
}
where $0<\eta \leq 1$.
\vspace{0.5cm}\\
Taking $\eta=\sqrt{\frac{2}{7}}$ and $\gamma$ sufficiently small, we get
\eq{\label{E14}
\arr{l}{\displaystyle\partial_\tau E_1 + 3 E_1 +\frac{1}{4}\left\|\nabla R\right\|^2_{L^2}\leq \left(\frac{35}{8}+CM\gamma \left(\frac{3}{2}-\theta\right)^2\right)\left\|\left(-\Delta\right)^{-\left(\frac{1+\theta}{2}\right)}R\right\|^2_{L^2}\\
\hspace{1cm}+ C M\gamma \left(\frac{3}{2}-\theta\right)^2\left(\left\|\left(-\Delta\right)^{-\left(\frac{1+\theta}{2}\right)} R\right\|^2_{L^2}+\alpha^2 e^{-2\tau}\left\|\Delta R\right\|^2_{L^2}\right)+ C  M^2 \gamma \left(\frac{3}{2}-\theta\right)^2 e^{-4\tau}.
}
}
Using the two energies $E_0$ and $E_1$, we define
\cent{$
\displaystyle E_2 = 6 E_0 + E_1.
$}
Combining the inequalities (\ref{eqE0}) and (\ref{E14}) and setting $\gamma$ sufficiently small, it is easy to check that
\eq{\label{eqE2}
\arr{l}{\displaystyle \partial E_2(\tau)+2\theta E_2(\tau)+\left\|\left(-\Delta\right)^{-\left(\frac{1+\theta}{2}\right)} R\right\|^2_{L^2} + \frac{1}{4}\left\|\nabla R\right\|^2_{L^2} \leq \\
\hspace{3cm} CM\gamma \left(\left\|\left|X\right|^4 R\right\|^2_{L^2}+\alpha^2 e^{-2\tau}\left\|\Delta R\right\|^2_{L^2}\right)+C M^2 \gamma \left(\frac{3}{2}-\theta\right) e^{-4\tau}.}
}
\subsection{Estimates in $H^2(\R^3)$}
In this part, we perform an $H^2-$estimate for the solution $R$ of (\ref{g2Re}) under the smallness assumption (\ref{cond2W}). To this end, we consider the $L^2-$scalar product of (\ref{g2Re}) with $-\Delta R$. We define the functional 
\cent{
$
E_3(\tau) = \frac{1}{2}\left(\left\|\nabla R\right\|^2_{L^2}+ \alpha e^{-\tau} \left\|\Delta R\right\|^2_{L^2}\right).
$
}
The next lemma gives the estimate of $R$ in the space $H^2(\R^3)$.
\lem{\label{lemE3}
Let $W\in  C^1\left(\left(\tau_0,\tau_\varepsilon\right), \mathbb H^1(4)\right)\cap C^0\left(\left(\tau_0,\tau_\varepsilon\right), \mathbb H^3(4)\right)$ be the solution of (\ref{g2We}). There exist two positive constants $\gamma_0$ and $T_0$ such that, if $T\geq T_0$ and $W$ satisfies the condition (\ref{cond2W}) for some positive constant $\gamma$ such that $\gamma \leq \gamma_0$, then there exists $C>0$ such that, for all $\tau \in \left[\tau_0,\tau^*_\varepsilon\right)$,
\eq{\label{eqE3}
\arr{l}{\displaystyle \partial_\tau E_3+3E_3+\frac{1}{2}\left\|\Delta R\right\|^2_{L^2}\leq  \frac{9}{4}\left\|\nabla R\right\|^2_{L^2} +CM \gamma \left(\frac{3}{2}-\theta\right)^2\left(\left\|R\right\|^2_{L^2}+\left\|\nabla R\right\|^2_{L^2}\right)\\
\displaystyle \hspace{11cm}+CM^2 \gamma \left(\frac{3}{2}-\theta\right)^2 e^{-\frac{7\tau}{2}}.
}
}
}
\textbf{Proof:} The proof of Lemma \ref{lemE3} is made through the $L^2-$scalar product of (\ref{g2Re}) with $-\Delta R$. First of all, we remark that 
\cent{
$
\displaystyle \curl\left(\left(W-\alpha e^{-\tau}\Delta W\right)\wedge U\right)=U.\nabla \left(W-\alpha e^{-\tau}\Delta W\right)-\left(W-\alpha e^{-\tau}\Delta W\right).\nabla U.
$
}
Making some computations that we let to the reader involving integrations by parts and the divergence free property of $U$, we obtain
\eq{\label{E31}
\partial_\tau\left(\left\|\nabla R\right\|^2_{L^2}+\alpha e^{-\tau}\left\|\Delta R\right\|^2_{L^2}\right)+\left(1-\frac{3\alpha}{4}e^{-\tau}\right)\left\|\Delta R\right\|^2_{L^2}=\frac{3}{4}\left\|\nabla R\right\|^2_{L^2}+I_1+I_2+I_3,
}
where
\cent{$
\arr{l}{
I_1 = \left(-U.\nabla \left(W-\alpha e^{-\tau}\Delta W\right), \Delta R\right)_{L^2},\\
\\
I_2= \left(\left(W-\alpha e^{-\tau}\Delta W\right).\nabla U, \Delta R\right)_{L^2},\\
\\
I_3= e^{-2\tau}\left(\alpha \Delta \Omega_\infty+\varepsilon\Delta^2 \Omega_\infty, \Delta R\right)_{L^2}.
}$
}
Like in the previous estimates, the easiest term is $I_3$. Indeed, using Hölder and Young inequalities and the inequality (\ref{cond2}), one has
\eq{\label{E3I3}
\arr{ll}{
\displaystyle I_3 &\displaystyle\leq e^{-2\tau}\left(\alpha \left\|\Delta \Omega_\infty\right\|_{L^2}+\varepsilon \left\|\Delta^2 \Omega_\infty\right\|_{L^2}\right) \left\|\Delta R\right\|_{L^2}\\
 & \displaystyle\leq \mu \left\|\Delta R\right\|^2 +\frac{C M \gamma \left(\frac{3}{2}-\theta\right)^2}{\mu} e^{-4\tau},
}
}
where $\mu$ is a positive constant which will be made more precise later.
\vspace{0.5cm}\\
We now look for an estimate of $I_1$. We decompose it as follows:
\cent{
$
I_1= I^1_1+I^2_1+I^3_1,
$
}
where
\cent{
$
\arr{l}{
\displaystyle I^1_1 =- e^{-\tau}\left(K.\nabla\left(\Omega_\infty-\alpha e^{-\tau}\Delta \Omega_\infty\right), \Delta R\right)_{L^2},\\
\\
\displaystyle I^2_1 = -e^{-2\tau}\left(V_\infty.\nabla\left(\Omega_\infty-\alpha e^{-\tau}\Delta \Omega_\infty\right), \Delta R\right)_{L^2},\\
\\
\displaystyle I^3_1 = -\left(U.\nabla\left(R-\alpha e^{-\tau}\Delta R\right), \Delta R\right)_{L^2}.
}
$
}
Due to the smoothness of $\Omega_\infty$ and the inequality (\ref{biotsb}), we get
\aligne{
I^1_1 &\leq e^{-\tau} \left\|K\right\|_{L^\infty} \left(\left\|\nabla \Omega_\infty\right\|_{L^2}+\alpha e^{-\tau}\left\|\nabla \Delta \Omega_\infty\right\|_{L^2}\right)\left\|\Delta R\right\|_{L^2}\\
&\leq C \left|b\right|e^{-\tau}\left\|R\right\|^{1/2}_{L^2}\left\|R\right\|^{1/2}_{L^6} \left\|\Delta R\right\|_{L^2}.
}
The continuous injection of $H^1(\R^3)$ into $L^6(\R^3)$, Young inequality and the inequality (\ref{cond2}) yield
\aligne{
I^1_1 &\leq C \left|b\right| e^{-\tau}\left\|R\right\|_{H^1} \left\|\Delta R\right\|_{L^2}\\
&\leq \mu \left\|\Delta R\right\|^2_{L^2} +\frac{CM\gamma\left(\frac{3}{2}-\theta\right)^2}{\mu}e^{-2\tau} \left(\left\|R\right\|^2_{L^2}+\left\|\nabla R\right\|^2_{L^2}\right).
}
Doing the same computations, we get
\aligne{
I^2_1 &\leq \mu \left\|\Delta R\right\|^2_{L^2}+\frac{CM^2\gamma^2\left(\frac{3}{2}-\theta\right)^4}{\mu}e^{-4\tau}.
}
The divergence free property of $U$ and an integration by parts imply
\cent{$
I^3_1 = \left(U.\nabla R, \Delta R\right)_{L^2}.
$}
Thus, using the Hölder and Young inequalities, Lemma \ref{biots} and the inequality (\ref{cond2W}), we obtain
\aligne{
I^3_1 & \leq \left\|U\right\|_{L^\infty} \left\|\nabla R\right\|_{L^2}\left\|\Delta R\right\|_{L^2}\\
&\leq C \left\|W\right\|^{1/2}_{L^2} \left\|W\right\|^{1/2}_{L^6}\left\|\nabla R\right\|_{L^2}\left\|\Delta R\right\|_{L^2}\\
&\leq C \left\|W\right\|_{H^1}\left\|\nabla R\right\|_{L^2}\left\|\Delta R\right\|_{L^2}\\
&\leq \mu \left\|\Delta R\right\|^2_{L^2}+\frac{CM\gamma \left(\frac{3}{2}-\theta\right)^2}{\mu}\left\|\nabla R\right\|^2_{L^2}.
}
Consequently, we have shown that
\eq{\label{E3I1}
\displaystyle I_1 \leq 3 \mu \left\|\Delta R\right\|^2_{L^2} +\frac{C M\gamma \left(\frac{3}{2}-\theta\right)^2}{\mu} \left(\left\|R\right\|^2_{L^2}+\left\|\nabla R\right\|^2_{L^2}\right)+ \frac{C M^2\gamma^2 \left(\frac{3}{2}-\theta\right)^4}{\mu}e^{-4\tau}.
}
It remains to estimate $I_2$. We set
\cent{
$
I_2 = I^1_2+I^2_2,
$}
where
\cent{
$
\arr{l}{
I^1_2 = -\left(W.\nabla U, \Delta R\right)_{L^2},\\
\\
I^2_2 = \alpha e^{-\tau}\left(\Delta W.\nabla U,\Delta R\right)_{L^2}.
}
$
}
Recalling that $W=e^{-\tau}\Omega_\infty+R$ and using Hölder and Young inequalities and the inequality (\ref{biotsc}) with $p=4$, one has
\aligne{
I^1_2 &\leq \left\|W\right\|_{L^4}\left\|\nabla U\right\|_{L^4}\left\|\Delta R\right\|_{L^2}\\
&\leq C\left\|W\right\|^2_{L^4}\left\|\Delta R\right\|_{L^2}\\
&\leq \mu \left\|\Delta R\right\|^2_{L^2} +\frac{C}{\mu}\left\|W\right\|^4_{L^4}\\
&\leq \mu \left\|\Delta R\right\|^2_{L^2} +\frac{C}{\mu}\left(e^{-4\tau} \left\|\Omega_\infty\right\|^4_{L^4}+\left\|R\right\|^4_{L^4}\right).
}
The condition (\ref{cond2}) and the continuous injection of $H^1(\R^3)$ into $L^4(\R^3)$ yield
\aligne{
I^1_2 &\leq \mu \left\|\Delta R\right\|^2_{L^2} +\frac{CM^2\gamma^2 \left(\frac{3}{2}-\theta\right)^4}{\mu}e^{-4\tau}+\frac{C}{\mu}\left\|R\right\|^4_{H^1}\\
&\leq \mu \left\|\Delta R\right\|^2_{L^2} +\frac{C M^2\gamma^2 \left(\frac{3}{2}-\theta\right)^4}{\mu}e^{-4\tau}+\frac{CM \gamma \left(\frac{3}{2}-\theta\right)^2}{\mu}\left(\left\|R\right\|^2_{L^2}+\left\|\nabla R\right\|^2_{L^2}\right).
}
Using the inequality (\ref{biotsb}) with $p=2$, $q=6$ and $\eta=\frac{1}{2}$ and the continuous injection of $H^1(\R^3)$ into $L^6(\R^3)$, we obtain
\aligne{
I^2_2 &\leq \alpha e^{-\tau}\left(\left\|\Delta R\right\|_{L^2}+e^{-\tau}\left\|\Delta \Omega_\infty \right\|_{L^2}\right)\left\|\nabla U\right\|_{L^\infty}\left\|\Delta R\right\|_{L^2}\\
&\leq C\alpha e^{-\tau}\left(\left\|\Delta R\right\|_{L^2}+e^{-\tau}\left\|\Delta \Omega_\infty \right\|_{L^2}\right)\left\|\nabla W\right\|^{1/2}_{L^2}\left\|\nabla W\right\|^{1/2}_{L^6}\left\|\Delta R\right\|_{L^2}\\
&\leq C\alpha e^{-\tau}\left(\left\|\Delta R\right\|_{L^2}+e^{-\tau}\left\|\Delta \Omega_\infty \right\|_{L^2}\right)\left\|\nabla W\right\|^{1/2}_{L^2}\left\|W\right\|^{1/2}_{H^2}\left\|\Delta R\right\|_{L^2}.
}
We set $\delta = M \gamma \left(\frac{3}{2}-\theta\right)^2$. Taking into account the inequalities (\ref{cond2}) and (\ref{cond2W}), it comes,
\aligne{
I^2_2 &\leq C \delta^{1/2} e^{-\frac{3\tau}{4}}\left(\left\|\Delta R\right\|_{L^2}+\delta^{1/2}e^{-\tau}\right) \left\|\Delta R\right\|_{L^2}\\
&\leq C \delta^{1/2}e^{-\frac{3\tau}{4}}\left\|\Delta R\right\|^2_{L^2}+C\delta e^{-\frac{7\tau}{4}}\left\|\Delta R\right\|_{L^2}\\
&\leq C \left(\delta^{1/2}e^{-\frac{3\tau}{4}}+\delta\right)\left\|\Delta R\right\|^2_{L^2}+C \delta e^{-\frac{7\tau}{2}}\\
&\leq C M \gamma^{1/2} \left(\frac{3}{2}-\theta\right) \left\|\Delta R\right\|^2_{L^2}+C M \gamma \left(\frac{3}{2}-\theta\right)^2 e^{-\frac{7\tau}{2}}.
}
Finally, we have shown,
\eq{
\label{E3I2}
\arr{l}{\displaystyle I_2 \leq \left(C M \gamma^{1/2} \left(\frac{3}{2}-\theta\right)+\mu\right) \left\|\Delta R\right\|^2_{L^2} +\frac{CM \gamma \left(\frac{3}{2}-\theta\right)^2 }{\mu}\left(\left\|R\right\|^2_{L^2}+\left\|\nabla R\right\|^2_{L^2}\right)\\
\hspace{11cm}\displaystyle +\frac{C M^2 \gamma \left(\frac{3}{2}-\theta\right)^2}{\mu}e^{-\frac{7\tau}{2}}.
}
}
Going back to (\ref{E31}), the inequalities (\ref{E3I3}), (\ref{E3I1}) and (\ref{E3I2}) imply
\eq{
\arr{l}{\displaystyle \partial_\tau E_3+3E_3+\left(1-5 \mu -\frac{9\alpha}{4}e^{-\tau}\right)\left\|\Delta R\right\|^2_{L^2}\leq \frac{9}{4}\left\|\nabla R\right\|^2_{L^2}+C M \gamma^{1/2} \left(\frac{3}{2}-\theta\right) \left\|\Delta R\right\|^2_{L^2}\\
\hspace{4.5cm}\displaystyle +\frac{CM \gamma \left(\frac{3}{2}-\theta\right)^2}{\mu}\left(\left\|R\right\|^2_{L^2}+\left\|\nabla R\right\|^2_{L^2}\right) +\frac{CM^2 \gamma \left(\frac{3}{2}-\theta\right)^2}{\mu}e^{-\frac{7\tau}{2}}.
}
}
We take $\gamma_0$ and $\mu$ small enough and $T=e^{\tau_0}$ large enough compared to $\alpha$ and obtain
\eq{
\arr{l}{\displaystyle \partial_\tau E_3+3E_3+\frac{1}{2}\left\|\Delta R\right\|^2_{L^2}\leq  \frac{9}{4}\left\|\nabla R\right\|^2_{L^2} +CM \gamma \left(\frac{3}{2}-\theta\right)^2\left(\left\|R\right\|^2_{L^2}+\left\|\nabla R\right\|^2_{L^2}\right)\\
\displaystyle \hspace{11cm}+CM^2 \gamma \left(\frac{3}{2}-\theta\right)^2 e^{-\frac{7\tau}{2}}.
}
}
\begin{flushright}
$\square$
\end{flushright}
To achieve the $H^2-$estimate, we combine $E_2$ and $E_3$ to define the functional
\cent{
$
E_4=12 E_2+E_3.
$}
Taking into account the two inequalities (\ref{eqE2}) and (\ref{eqE3}), we see that $E_4$ satisfies
\eq{
\arr{l}{\displaystyle\partial_\tau E_4+2\theta E_4 +12\left\|\left(-\Delta\right)^{-\left(\theta-\frac{1}{4}\right)} R\right\|^2_{L^2} +\frac{3}{4} \left\|\nabla R\right\|^2_{L^2}+ \frac{1}{2}\left\|\Delta R\right\|^2_{L^2} \leq \\
\\
\hspace{1cm}\displaystyle+CM \gamma \left(\left\|R\right\|^2_{L^2}+\left\|\nabla R\right\|^2_{L^2}+\left\|\left|X\right|^4 R\right\|^2_{L^2}+\alpha^2 e^{-2\tau}\left\|\Delta R\right\|^2_{L^2}\right)+C M^2 \gamma \left(\frac{3}{2}-\theta\right) e^{-\frac{7\tau}{2}}.
}
}
Using again the interpolation inequality (\ref{interp}) and taking $\gamma_0$ small enough, this inequality becomes
\eq{\label{eqE4}
\arr{l}{\partial_\tau E_4+2\theta E_4 + 10\left\|\left(-\Delta\right)^{-\left(\theta-\frac{1}{4}\right)} R\right\|^2_{L^2} + \frac{1}{2}\left\|\nabla R\right\|^2_{L^2}+ \frac{1}{4}\left\|\Delta R\right\|^2_{L^2} \leq \\
\hspace{7cm} \displaystyle CM \gamma\left\|\left|X\right|^4 R\right\|^2_{L^2}+CM^2 \gamma \left(\frac{3}{2}-\theta\right) e^{-\frac{7\tau}{2}}.
}
}
\subsection{Estimates in $H^2(4)$}
To finish the energy estimates, we have to work in weighted spaces. We can see that the terms of the right hand side of the inequality (\ref{eqE4}) involve weighted $L^2-$norms that we have to absorb. In order to perform estimates in weighted Lebesgue norms, and additionally absorb the weighted terms of (\ref{eqE4}), we make the $L^2-$inner product of (\ref{g2Re}) with $\left|X\right|^8 \left(R-\alpha e^{-\tau} \Delta R\right)$. One defines the energy functional
\cent{
$
\displaystyle E_5 = \frac{1}{2}\left\|\left|X\right|^4 \left(R-\alpha e^{-\tau}\Delta R\right)\right\|^2_{L^2}.
$}
The next lemma summarizes the terms provided by the linear part of (\ref{g2Re}), when making the $L^2-$scalar product with $\left|X\right|^8 \left(R-\alpha e^{-\tau} \Delta R\right)$.
\lem{\label{lemcalc}
Let $u$ be a divergence free vector field of $\mathbb H^2(4)$, $a\in \R$ and $F(u)=\left|x\right|^8 \left(u- a \Delta u\right)$. The five next equalities hold.
\begin{enumerate}
\item \eq{\label{calc1}\displaystyle \left(\Delta u, F(u)\right)_{L^2}=36\left\|\left|x\right|^3 u\right\|^2_{L^2}-\left\|\left|x\right|^4 \nabla u\right\|^2_{L^2}-a\left\|\left|x\right|^4 \Delta u\right\|^2_{L^2}.
}
\item  \eq{\label{calc2}\displaystyle \left(\frac{x}{2}.\nabla u, F(u)\right)_{L^2}=-\frac{11}{4}\left\|\left|x\right|^4 u\right\|^2_{L^2}-\frac{9 a}{4}\left\|\left|x\right|^4\nabla u\right\|^2_{L^2}+4a\left\|\left|x\right|^3\left(x.\nabla u\right)\right\|^2_{L^2}.
}
\item \eq{\label{calc3}\arr{l}{\displaystyle \left(\mathcal L (u), F(u)\right)_{L^2} = -\frac{7}{4}\left\|\left|x\right|^4 u\right\|^2_{L^2}-\left(1+\frac{5a}{4}\right)\left\|\left|x\right|^4\nabla u\right\|^2_{L^2}-a\left\|\left|x\right|^4\Delta u\right\|^2_{L^2}\\
\hspace{6cm}\displaystyle+4a\left\|\left|x\right|^3\left(x.\nabla u\right)\right\|^2_{L^2}+36\left(1-a\right) \left\|\left|x\right|^3 u\right\|^2_{L^2}.
}}
\item \eq{\label{calc4} \arr{ll}{\displaystyle \left(\Delta^2 u,F(u)\right)_{L^2}&=\left\|\left|x\right|^4 \Delta u\right\|^2_{L^2}-16\left\|\left|x\right|^3\nabla u\right\|^2_{L^2}-96\left\|\left|x\right|^2 \left(x.\nabla u\right)\right\|^2_{L^2}\\
&\hspace{1cm}\displaystyle +1512 \left\|\left|x\right|^2 u\right\|^2_{L^2}+a\left\|\left|x\right|^4 \nabla \Delta u\right\|^2_{L^2}-36a\left\|\left|x\right|^3 \Delta u\right\|^2_{L^2}.}
}
\item \eq{\label{calc5}\arr{l}{\displaystyle\left(\frac{x}{2}.\nabla \Delta u,F(u)\right)_{L^2}=\frac{13}{4}\left\|\left|x\right|^4 \nabla u\right\|^2_{L^2}+\frac{11a}{4} \left\|\left|x\right|^4  \Delta u\right\|^2_{L^2}\\
\hspace{6cm}+4\left\|\left|x\right|^3\left(x.\nabla u\right)\right\|^2_{L^2}-180  \left\|\left|x\right|^3 u\right\|^2_{L^2}.}}
\end{enumerate}
}
There is no difficulty in the proof of this lemma, which is let to the reader. It is only a consequence of many integrations by parts.
\vspace{0.5cm}\\
The next lemma enables us to close the $\mathbb H^2(4)$ estimate.
\lem{\label{lemE5}
Let $W\in C^1\left(\left(\tau_0,\tau_\varepsilon\right),\mathbb H^1(4)\right)\cap C^0\left(\left(\tau_0,\tau_\varepsilon\right),\mathbb H^3(4)\right)$ be the solution of (\ref{g2We}). There exist two positive constants $\gamma_0$ and $T_0$ such that, if $T\geq T_0$ and $W$ satisfy the condition (\ref{cond2W}) for some positive constant  such that $\gamma \leq \gamma_0$, then there exists $C>0$ such that, for all $\tau \in \left[\tau_0,\tau^*_\varepsilon\right)$,
\eq{\label{eqE5}\arr{l}{
\displaystyle \partial_\tau E_5+3E_5+\frac{1}{16}\left\|\left|X\right|^4 R\right\|^2_{L^2}+\left(\frac{\alpha}{2} e^{-\tau}+\frac{\alpha^2}{4} e^{-2\tau}\right)\left\|\left|X\right|^4 \Delta R\right\|^2_{L^2} \leq  K_1\left\| R\right\|^2_{L^2}\\
\hspace{1cm} \displaystyle + CM \gamma^{1/4} \left(\frac{3}{2}-\theta\right)^{1/2}\left(\left\| R\right\|^2_{L^2}+\left\| \nabla R\right\|^2_{L^2}+\left\| \Delta R\right\|^2_{L^2}\right)+C M^2 \gamma \left(\frac{3}{2}-\theta\right)^{2} e^{-4\tau},
}
}
where $K_1$ is a positive constant independent of the parameters.
}
\textbf{Proof:} To obtain the inequality (\ref{eqE5}) of this lemma, we perform the $L^2-$inner product of (\ref{g2Re}) with $\left|X\right|^8\left(R-\alpha e^{-\tau}\Delta R\right)$. We deliberately omit the positive terms obtained from $\varepsilon \Delta^2 W$ which do not play any role in the next estimates. Using Lemma (\ref{lemcalc}) and making some easy computations, one obtains
\eq{\label{eqE51}\arr{l}{
\displaystyle \frac{1}{2}\partial_\tau\left(\left\|\left|X\right|^4 \left(R-\alpha e^{-\tau} \Delta R\right)\right\|^2_{L^2}\right)+\frac{7}{4}\left\|\left|X\right|^4 R\right\|^2_{L^2}+\left(1+\frac{7\alpha}{2}e^{-\tau}\right)\left\|\left|X\right|^4\nabla R\right\|^2_{L^2}\\
\\
\displaystyle \hspace{3.5cm}+\left(\alpha e^{-\tau}+\frac{7\alpha^2}{4} e^{-2\tau}\right)\left\|\left|X\right|^4 \Delta R\right\|^2_{L^2}-108 \alpha e^{-\tau}\left\|\left|X\right|^3 R\right\|^2_{L^2}=\\
\\
\hspace{9cm}36\left\|\left|X\right|^3 R\right\|^2_{L^2}+I_1+I_2+I_3+I_4,
}
}
where
\cent{
$
\arr{l}{
I_1=\left(-U.\nabla\left(W-\alpha e^{-\tau}\Delta W\right),\left|X\right|^8\left(R-\alpha e^{-\tau}\Delta R\right)\right)_{L^2},\\
\\
I_2=\left(\left(W-\alpha e^{-\tau}\Delta W\right).\nabla U,\left|X\right|^8\left(R-\alpha e^{-\tau}\Delta R\right)\right)_{L^2},\\
\\
I_3=\left(-\varepsilon e^{-2\tau} \Delta^2 \Omega_\infty-\alpha e^{-2\tau}\Delta \Omega_\infty,\left|X\right|^8\left(R-\alpha e^{-\tau}\Delta R\right)\right)_{L^2},\\
\\
I_4 = \varepsilon e^{-\tau} \left(16\left\|\left|X\right|^3\nabla R\right\|^2_{L^2}+96\left\|\left|X\right|^2 \left(X.\nabla R\right)\right\|^2_{L^2}+36\alpha e^{-\tau}\left\|\left|X\right|^3 \Delta R\right\|^2_{L^2}\right).
}
$
}
In the proof of this lemma, we use the notation
\cent{
$
\delta = M \gamma \left(\frac{3}{2}-\theta\right)^2.
$
}
As usual, $I_3$ is the easiest term to estimate. Indeed, due to the smoothness of $\Omega_\infty$ and the inequality (\ref{cond2}), we get
\eq{\label{eqE5I3}\arr{ll}{
\displaystyle I_3 & \displaystyle\leq Ce^{-2\tau}\left\|\left|X\right|^4\left(\alpha \Delta \Omega_\infty +\varepsilon \Delta^2 \Omega_\infty\right)\right\|_{L^2}\left(\left\|\left|X\right|^4 R\right\|_{L^2}+\alpha e^{-\tau}\left\|\left|X\right|^4\Delta R\right\|_{L^2}\right)\\
&\displaystyle \leq \mu\left\|\left|X\right|^4 R\right\|^2_{L^2}+\mu\alpha^2 e^{-2\tau}\left\|\left|X\right|^4\Delta R\right\|^2_{L^2} + \frac{C\left|b\right|^2}{\mu} e^{-4\tau}\\
&\displaystyle \leq \mu\left\|\left|X\right|^4 R\right\|^2_{L^2}+\mu\alpha^2 e^{-2\tau}\left\|\left|X\right|^4\Delta R\right\|^2_{L^2} + \frac{CM \gamma \left(\frac{3}{2}-\theta\right)^2}{\mu} e^{-4\tau},
}
}
where $\mu$ is a positive constant that will be made more precise later.
\vspace{0.5cm}\\
We now give an estimate of $I_4$, which is also quite simple to bound. We just need Hölder and Young inequalities to estimate this term in a convenient way. Indeed, using convexity inequalities, it is simple to show that
\cent{$\left\|\left|X\right|^3\nabla R\right\|^2_{L^2} +\left\|\left|X\right|^2 \left(X.\nabla R\right)\right\|^2_{L^2}\leq C \left\|\left|X\right|^4 \nabla R\right\|^2_{L^2}+C\left\|\nabla R\right\|^2_{L^2},$}
and
\cent{$\alpha e^{-\tau}\left\|\left|X\right|^3 \Delta R\right\|^2_{L^2} \leq C \alpha e^{-\tau}\left\|\left|X\right|^4 \Delta R\right\|^2_{L^2}+C\alpha e^{-\tau}\left\|\Delta R\right\|^2_{L^2}.$
}
Thus, if we take $\varepsilon \leq \alpha M\gamma\left(\frac{3}{2}-\theta\right)^2$, we get
\eq{\label{eqE5I4}
\arr{l}{I_4 \leq CM\gamma\left(\frac{3}{2}-\theta\right)^2\left(\alpha e^{-\tau}\left\|\left|X\right|^4 \nabla R\right\|^2_{L^2} + \alpha^2 e^{-2\tau}\left\|\left|X\right|^4 \Delta R\right\|^2_{L^2}\right)\\
\hspace{5cm}+CM\gamma\left(\frac{3}{2}-\theta\right)^2\left(\alpha e^{-\tau}\left\|\nabla R\right\|^2_{L^2} + \alpha^2 e^{-2\tau}\left\|\Delta R\right\|^2_{L^2}\right).}
}
As for the $H^2-$estimate, we have to study separately $I_1$ and $I_2$. We begin with $I_1$, that we rewrite
\cent{
$
I_1 = I^1_1+I^2_1+I^3_1,
$
}
where
\cent{
$\arr{l}{
\displaystyle I^1_1 = \left(U.\nabla\left(R-\alpha e^{-\tau}\Delta R\right),\left|X\right|^8\left(R-\alpha e^{-\tau}\Delta R\right)\right)_{L^2},\\
\\
\displaystyle I^2_1 = e^{-2\tau}\left(V_\infty.\nabla\left(\Omega_\infty-\alpha e^{-\tau}\Delta \Omega_\infty\right),\left|X\right|^8\left(R-\alpha e^{-\tau}\Delta R\right)\right)_{L^2},\\
\\
\displaystyle I^3_1 = e^{-\tau}\left(K.\nabla\left(\Omega_\infty-\alpha e^{-\tau}\Delta \Omega_\infty\right),\left|X\right|^8\left(R-\alpha e^{-\tau}\Delta R\right)\right)_{L^2}.
}$
}
Using an integration by parts, the fact that $\dv U = 0$ and the Hölder inequalities, one has
\aligne{
I^1_1 & = \frac{1}{2}\int_{\R^3} \left|X\right|^8 U(X).\nabla\left(\left|R(X)-\alpha e^{-\tau}\Delta R(X)\right|^2\right) dX\\
& = -4\int_{\R^3} \left|X\right|^6 \left(X. U(X)\right)\left|R(X)-\alpha e^{-\tau}\Delta R(X)\right|^2 dX\\
&\leq C\left\|U\right\|_{L^\infty}\left(\left\|\left|X\right|^{7/2} R\right\|^2_{L^2}+\alpha^2 e^{-2\tau} \left\|\left|X\right|^{7/2} \Delta R\right\|^2_{L^2}\right).
}
The inequalities (\ref{biotsb}) with $p=2$, $q=6$ and $\eta = \frac{1}{2}$ and (\ref{cond2W}) and the continuous injection of $H^1(\R^3)$ into $L^6(\R^3)$ imply
\aligne{
I^1_1&\leq C\left\|W\right\|^{1/2}_{L^2}\left\|W\right\|^{1/2}_{L^6}\left(\left\|R\right\|^2_{L^2}+\left\|\left|X\right|^4 R\right\|^2_{L^2}+\alpha^2 e^{-2\tau} \left\|\Delta R\right\|^2_{L^2}+\alpha^2 e^{-2\tau} \left\|\left|X\right|^4 \Delta R\right\|^2_{L^2}\right)\\
&\leq C\left\|W\right\|_{H^1}\left(\left\|R\right\|^2_{L^2}+\left\|\left|X\right|^4 R\right\|^2_{L^2}+\alpha^2 e^{-2\tau} \left\|\Delta R\right\|^2_{L^2}+\alpha^2 e^{-2\tau} \left\|\left|X\right|^4 \Delta R\right\|^2_{L^2}\right)\\
&\leq CM^{1/2}\gamma^{1/2}\left(\frac{3}{2}-\theta\right)\Big(\left\|R\right\|^2_{L^2}+\left\|\left|X\right|^4 R\right\|^2_{L^2}\\
&\hspace{4cm}+\alpha^2 e^{-2\tau} \left\|\Delta R\right\|^2_{L^2}+\alpha^2 e^{-2\tau} \left\|\left|X\right|^4 \Delta R\right\|^2_{L^2}\Big).\\
}
Because of the smoothness of $\Omega_\infty$, $I^2_1$ is a little easier to estimate. Indeed using once more the inequalities (\ref{biotsb}) and (\ref{cond2}) and the Hölder and Young inequalities, we get
\aligne{
I^2_1 & \leq C e^{-2\tau}\left\|V_\infty\right\|_{L^\infty} \left(\left\|\left|X\right|^4 \nabla \Omega_\infty \right\|_{L^2}+\alpha e^{-\tau}\left\|\left|X\right|^4 \nabla \Delta \Omega_\infty\right\|_{L^2}\right) \\
&\hspace{6cm}\left(\left\|\left|X\right|^4 R\right\|_{L^2}+\alpha e^{-\tau}\left\|\left|X\right|^4 \Delta R\right\|_{L^2}\right)\\
& \leq C \left|b\right|e^{-2\tau}\left\|\Omega_\infty\right\|^{1/2}_{L^2} \left\|\Omega_\infty\right\|^{1/2}_{L^6}\left(\left\|\left|X\right|^4 R\right\|_{L^2}+\alpha e^{-\tau}\left\|\left|X\right|^4 \Delta R\right\|_{L^2}\right)\\
&\leq C \left|b\right|^2e^{-2\tau}\left(\left\|\left|X\right|^4 R\right\|_{L^2}+\alpha e^{-\tau}\left\|\left|X\right|^4 \Delta R\right\|_{L^2}\right)\\
&\leq \mu\left\|\left|X\right|^4 R\right\|^2_{L^2}+\mu \alpha^2 e^{-2\tau}\left\|\left|X\right|^4 \Delta R\right\|^2_{L^2}+\frac{CM^2\gamma^2\left(\frac{3}{2}-\theta\right)^4}{\mu} e^{-4\tau}.\\
}
Likewise, we get
\aligne{
I^3_1 & \leq C\left|b\right| e^{-\tau}\left\|K\right\|_{L^\infty}\left(\left\|\left|X\right|^4 R\right\|_{L^2}+\alpha e^{-\tau}\left\|\left|X\right|^4 \Delta R\right\|_{L^2}\right)\\
&\leq C \left|b\right| e^{-\tau}\left\|R\right\|_{H^1}\left(\left\|\left|X\right|^4 R\right\|_{L^2}+\alpha e^{-\tau}\left\|\left|X\right|^4 \Delta R\right\|_{L^2}\right)\\
&\leq C M^{1/2}\gamma^{1/2}\left(\frac{3}{2}-\theta\right)\left(\left\|R\right\|^2_{L^2}+\left\|\nabla R\right\|^2_{L^2}+\left\|\left|X\right|^4 R\right\|^2_{L^2}+\alpha^2 e^{-2\tau}\left\|\left|X\right|^4 \Delta R\right\|^2_{L^2}\right).\\
}
Finally, taking $T$ so that $\displaystyle \alpha e^{-\tau_0}= \frac{\alpha}{T}\leq 1$, we have
\eq{\label{eqE5I1}
\arr{l}{\displaystyle I_1 \leq \mu\left\|\left|X\right|^4 R\right\|^2_{L^2}+\mu \alpha e^{-\tau}\left\|\left|X\right|^4 \Delta R\right\|^2_{L^2}+\frac{CM^2\gamma^2\left(\frac{3}{2}-\theta\right)^4}{\mu} e^{-4\tau}\\
\displaystyle \hspace{2cm}+ CM^{1/2}\gamma^{1/2}\left(\frac{3}{2}-\theta\right)\Big(\left\|R\right\|^2_{L^2}+\left\|\nabla R\right\|^2_{L^2}+ \left\|\Delta R\right\|^2_{L^2}\\
\displaystyle \hspace{7cm}+\left\|\left|X\right|^4 R\right\|^2_{L^2}+\alpha e^{-\tau} \left\|\left|X\right|^4 \Delta R\right\|^2_{L^2}\Big)}
}
It remains to bound $I_2$, which is the hardest term to estimate. Like for $I_1$, we rewrite it
\cent{
$
I_2 = I^1_2+I^2_2+I^3_2+I^4_2,
$}
where
\cent{
$
\arr{l}{
I^1_2 = e^{-\tau}\left(\left(R-\alpha e^{-\tau}\Delta R\right).\nabla V_\infty,\left|X\right|^8\left(R-\alpha e^{-\tau}\Delta R\right)\right)_{L^2},\\
\\
I^2_2 = \left(\left(R-\alpha e^{-\tau}\Delta R\right).\nabla K,\left|X\right|^8\left(R-\alpha e^{-\tau}\Delta R\right)\right)_{L^2},\\
\\
I^3_2 = e^{-2\tau}\left(\left(\Omega_\infty-\alpha e^{-\tau}\Delta \Omega_\infty\right).\nabla V_\infty,\left|X\right|^8\left(R-\alpha e^{-\tau}\Delta R\right)\right)_{L^2},\\
\\
I^4_2 =  e^{-\tau}\left(\left(\Omega_\infty-\alpha e^{-\tau}\Delta \Omega_\infty\right).\nabla K,\left|X\right|^8\left(R-\alpha e^{-\tau}\Delta R\right)\right)_{L^2}.\\
}
$
}
Using the inequality (\ref{biotsb}) and the smoothness of $\Omega_\infty$, we get
\aligne{
I^1_2 &\leq Ce^{-\tau} \left\|\nabla V_\infty\right\|_{L^\infty}\left(\left\|\left|X\right|^4 R\right\|^2_{L^2}+\alpha^2 e^{-2\tau}\left\|\left|X\right|^4 \Delta R\right\|^2_{L^2}\right)\\
&\leq Ce^{-\tau} \left\|\nabla \Omega_\infty\right\|^{1/2}_{L^2}\left\|\nabla \Omega_\infty\right\|^{1/2}_{L^6}\left(\left\|\left|X\right|^4 R\right\|^2_{L^2}+\alpha^2 e^{-2\tau}\left\|\left|X\right|^4 \Delta R\right\|^2_{L^2}\right)\\
&\leq C \left|b\right|\left(\left\|\left|X\right|^4 R\right\|^2_{L^2}+\alpha^2 e^{-2\tau}\left\|\left|X\right|^4 \Delta R\right\|_{L^2}^2\right)\\
&\leq C M^{1/2} \gamma^{1/2} \left(\frac{3}{2}-\theta\right)\left(\left\|\left|X\right|^4 R\right\|^2_{L^2}+\alpha^2 e^{-2\tau}\left\|\left|X\right|^4 \Delta R\right\|_{L^2}^2\right).
}
We now estimate $I^2_2$. We recall the notation $\delta = M\gamma \left(\frac{3}{2}-\theta\right)^2$. Using again the inequality (\ref{biotsb}), the inequality (\ref{cond2}) and the continuous injection of $H^1(\R^3)$ into $L^6(\R^3)$, one has
\aligne{
I^2_2 &\leq \left\|\nabla K\right\|_{L^\infty}  \left(\left\|\left|X\right|^4 R\right\|^2_{L^2}+\alpha^2 e^{-2\tau}\left\|\left|X\right|^4 \Delta R\right\|^2_{L^2}\right)\\
&\leq C\left\|\nabla R\right\|^{1/2}_{L_2}\left\|\nabla R\right\|^{1/2}_{L_6}  \left(\left\|\left|X\right|^4 R\right\|^2_{L^2}+\alpha^2 e^{-2\tau}\left\|\left|X\right|^4 \Delta R\right\|^2_{L^2}\right)\\
&\leq C\delta^{1/4}\left\|\nabla R\right\|^{1/2}_{H^1}  \left(\left\|\left|X\right|^4 R\right\|^2_{L^2}+\alpha^2 e^{-2\tau}\left\|\left|X\right|^4 \Delta R\right\|^2_{L^2}\right)\\
&\leq C\delta^{1/2} \left(\left\|\left|X\right|^4 R\right\|^2_{L^2}+\alpha^2 e^{-2\tau}\left\|\left|X\right|^4 \Delta R\right\|^2_{L^2}\right)\\
&\hspace{3cm}+ C\delta^{1/4}\left\|\Delta R\right\|^{1/2}_{L^2}  \left(\left\|\left|X\right|^4 R\right\|^2_{L^2}+\alpha^2 e^{-2\tau}\left\|\left|X\right|^4 \Delta R\right\|^2_{L^2}\right)\\
&\leq C\delta^{1/2} \left(\left\|\left|X\right|^4 R\right\|^2_{L^2}+\alpha^2 e^{-2\tau}\left\|\left|X\right|^4 \Delta R\right\|^2_{L^2}\right)+ C\delta^{1/4}\left\|\Delta R\right\|^{1/2}_{L^2}  \left\|\left|X\right|^4 R\right\|^2_{L^2}\\
&\hspace{9cm}+C\delta^{1/2}\alpha^{7/4} e^{-\frac{7\tau}{4}}\left\|\left|X\right|^4 \Delta R\right\|^2_{L^2}.
}
To finish the estimate of $I^2_2$, we use the convexity inequality $ab \leq \frac{3}{4} a^{\frac{4}{3}}+\frac{1}{4}b^{4}$ and the condition (\ref{cond2}). We obtain
\aligne{
I^2_2 &\leq C\delta^{1/2}  \left(\left\|\left|X\right|^4 R\right\|^2_{L^2}+\alpha^2 e^{-2\tau}\left\|\left|X\right|^4 \Delta R\right\|^2_{L^2}\right)+ C\delta^{1/4}\left(\left\|\Delta R\right\|^2_{L^2}  +\left\|\left|X\right|^4 R\right\|^{8/3}_{L^2}\right)\\
&\hspace{9cm}+C\delta^{1/2}\alpha^{7/4} e^{-\frac{7\tau}{4}}\left\|\left|X\right|^4 \Delta R\right\|^2_{L^2}\\
&\leq C\delta^{1/2}  \left(\left\|\left|X\right|^4 R\right\|^2_{L^2}+\alpha^2 e^{-2\tau}\left\|\left|X\right|^4 \Delta R\right\|^2_{L^2}\right)+ C\delta^{1/4}\left\|\Delta R\right\|^2_{L^2}  +C\delta^{7/12}\left\|\left|X\right|^4 R\right\|^2_{L^2}\\
&\hspace{9cm}+C\delta^{1/2}\alpha^{7/4} e^{-\frac{7\tau}{4}}\left\|\left|X\right|^4 \Delta R\right\|^2_{L^2}.
}
Consequently, if we assume $\gamma \leq 1$ and $\left(\frac{3}{2}-\theta\right) \leq 1$, one has
\cent{
$
\displaystyle I^2_2  \leq C M^{7/12} \gamma^{1/4} \left(\frac{3}{2}-\theta\right)^{1/2} \left(\left\|\left|X\right|^4 R\right\|^2_{L^2}+\left\|\Delta R\right\|^2_{L^2}+\alpha^2 e^{-2\tau}\left\|\left|X\right|^4 \Delta R\right\|^2_{L^2} \right)
$
}
It it easier bound $I^3_2$. Indeed, the inequality (\ref{biotsc}) and the inequality (\ref{cond2}) imply
\aligne{
I^3_2 & \leq C e^{-2\tau}\left\|\Omega_\infty -\alpha e^{-\tau}\Delta \Omega_\infty\right\|_{L^\infty} \left\|\nabla V_\infty\right\|_{L^2}\left(\left\|\left|X\right|^4 R\right\|_{L^2}+\alpha e^{-\tau} \left\|\left|X\right|^4 \Delta R\right\|_{L^2}\right)\\
& \leq C \left|b\right|e^{-2\tau}\left\|\Omega_\infty\right\|_{L^2}\left(\left\|\left|X\right|^4 R\right\|_{L^2}+\alpha e^{-\tau} \left\|\left|X\right|^4 \Delta R\right\|_{L^2}\right)\\
& \leq C \left|b\right|^2 e^{-2\tau}\left(\left\|\left|X\right|^4 R\right\|_{L^2}+\alpha e^{-\tau} \left\|\left|X\right|^4 \Delta R\right\|_{L^2}\right)\\
& \leq \mu\left\|\left|X\right|^4 R\right\|^2_{L^2}+\mu\alpha^2 e^{-2\tau} \left\|\left|X\right|^4 \Delta R\right\|^2_{L^2}+\frac{CM^2\gamma^2\left(\frac{3}{2}-1\right)^4}{\mu} e^{-4\tau}.
}
Likewise, we obtain
\aligne{
I^4_2 & \leq C e^{-\tau} \left\|\Omega_\infty-\alpha e^{-\tau} \Delta \Omega_\infty\right\|_{L^\infty} \left\|\nabla K\right\|_{L^2}\left(\left\|\left|X\right|^4 R\right\|_{L^2}+\alpha e^{-\tau} \left\|\left|X\right|^4 \Delta R\right\|_{L^2}\right)\\
& \leq C \left|b\right|e^{-\tau} \left\| R\right\|_{L^2}\left(\left\|\left|X\right|^4 R\right\|_{L^2}+\alpha e^{-\tau} \left\|\left|X\right|^4 \Delta R\right\|_{L^2}\right)\\
& \leq C M^{1/2}\gamma^{1/2}\left(\frac{3}{2}-1\right) \left(\left\| R\right\|^2_{L^2}+\left\|\left|X\right|^4 R\right\|^2_{L^2}+\alpha^2 e^{-2\tau} \left\|\left|X\right|^4 \Delta R\right\|^2_{L^2}\right).\\
}
Thus, taking $T_0$ large enough so that $\alpha e^{-\tau_0}=\frac{\alpha}{T}\leq 1$, the following inequality holds:
\eq{\label{eqE5I2}
\arr{l}{
\displaystyle I_2 \leq \mu\left(\left\|\left|X\right|^4 R\right\|^2_{L^2}+\alpha e^{-\tau} \left\|\left|X\right|^4 \Delta R\right\|^2_{L^2}\right)+\frac{C M^2\gamma^2\left(\frac{3}{2}-1\right)^4}{\mu} e^{-4\tau}\\
\displaystyle\hspace{1cm}+CM\gamma^{1/4}\left(\frac{3}{2}-1\right)^{1/2}\left(\left\| R\right\|^2_{L^2}+\left\| \Delta R\right\|^2_{L^2}+\left\|\left|X\right|^4 R\right\|^2_{L^2}+\alpha e^{-\tau}\left\|\left|X\right|^4 \Delta R\right\|^2_{L^2}\right).
}
}
Combining the equality (\ref{eqE51}) together with the inequalities (\ref{eqE5I3}), (\ref{eqE5I4}), (\ref{eqE5I1}) and (\ref{eqE5I2}) and taking $T_0$ big enough compared to $\alpha$, we have
\eq{\label{eqE52}\arr{l}{
\displaystyle \frac{1}{2}\partial_\tau\left(\left\|\left|X\right|^4 \left(R-\alpha e^{-\tau} \Delta R\right)\right\|^2_{L^2}\right)+\frac{7}{4}\left\|\left|X\right|^4 R\right\|^2_{L^2}+\left(1+\frac{7\alpha}{2}e^{-\tau}\right)\left\|\left|X\right|^4\nabla R\right\|^2_{L^2}\\
\\
\hspace{4cm}\displaystyle+\left(\alpha e^{-\tau}+\frac{7\alpha^2}{4} e^{-2\tau}\right)\left\|\left|X\right|^4 \Delta R\right\|^2_{L^2}-108 \alpha e^{-\tau}\left\|\left|X\right|^3 R\right\|^2_{L^2} \leq\\
\\
\hspace{0cm}\displaystyle  C\left(M \gamma^{1/4} \left(\frac{3}{2}-\theta\right)^{1/2} +\mu\right)\left(\left\|\left|X\right|^4 R\right\|^2_{L^2}+\alpha e^{-\tau} \left\|\left|X\right|^4 \nabla R\right\|^2_{L^2}+\alpha e^{-\tau} \left\|\left|X\right|^4 \Delta R\right\|^2_{L^2}\right)\\
\\
\hspace{1cm}\displaystyle+36\left\|\left|X\right|^3 R\right\|^2_{L^2} +CM \gamma^{1/4} \left(\frac{3}{2}-\theta\right)^{1/2}\left(\left\| R\right\|^2_{L^2}+\left\|\nabla R\right\|^2_{L^2}+\left\| \Delta R\right\|^2_{L^2}\right)\\
\hspace{11cm}\displaystyle+\frac{C M^2 \gamma \left(\frac{3}{2}-\theta\right)^{2}}{\mu} e^{-4\tau}.\\
}
}
Integrating several times by parts, it is easy to check that
\eq{E_4=\frac{1}{2}\left\|\left|X\right|^4 R\right\|^2_{L^2}+\frac{\alpha^2}{2}e^{-2\tau}\left\|\left|X\right|^4 \Delta R\right\|^2_{L^2}+\alpha e^{-\tau}\left\|\left|X\right|^4 \nabla R\right\|^2_{L^2}-36 \alpha e^{-\tau}\left\|\left|X\right|^3 R\right\|^2_{L^2}.}
Consequently, the inequality (\ref{eqE52}) becomes
\eq{\label{eqE53}\arr{l}{
\displaystyle \partial_\tau E_5+3E_5+\frac{1}{4}\left\|\left|X\right|^4 R\right\|^2_{L^2}+\left(1+\frac{\alpha}{2}e^{-\tau}\right)\left\|\left|X\right|^4\nabla R\right\|^2_{L^2}\\
\\
\hspace{7cm}\displaystyle +\left(\alpha e^{-\tau}+\frac{\alpha^2}{4} e^{-2\tau}\right)\left\|\left|X\right|^4 \Delta R\right\|^2_{L^2} \leq \\
\\
\hspace{0cm}\displaystyle  C\left(M \gamma^{1/4} \left(\frac{3}{2}-\theta\right)^{1/2} +\mu\right)\left(\left\|\left|X\right|^4 R\right\|^2_{L^2}+\alpha e^{-\tau} \left\|\left|X\right|^4 \nabla R\right\|^2_{L^2}+\alpha e^{-\tau} \left\|\left|X\right|^4 \Delta R\right\|^2_{L^2}\right)\\
\\
\hspace{1cm}\displaystyle+36\left\|\left|X\right|^3 R\right\|^2_{L^2} +CM \gamma^{1/4} \left(\frac{3}{2}-\theta\right)^{1/2}\left(\left\| R\right\|^2_{L^2}+\left\|\nabla R\right\|^2_{L^2}+\left\| \Delta R\right\|^2_{L^2}\right)\\
\hspace{11cm}\displaystyle+\frac{C M^2 \gamma \left(\frac{3}{2}-\theta\right)^{2}}{\mu} e^{-4\tau}.\\
}
}
Thus, taking $\gamma_0$ and $\mu$ small enough, we obtain
\eq{\label{eqE54}\arr{l}{
\displaystyle \partial_\tau E_5+3E_5+\frac{1}{8}\left\|\left|X\right|^4 R\right\|^2_{L^2}+\left(\frac{\alpha}{2} e^{-\tau}+\frac{\alpha^2}{4} e^{-2\tau}\right)\left\|\left|X\right|^4 \Delta R\right\|^2_{L^2} \leq 36\left\|\left|X\right|^3 R\right\|^2_{L^2}\\
\\
\hspace{1.5cm}  +\displaystyle CM \gamma^{1/4} \left(\frac{3}{2}-\theta\right)^{1/2}\left(\left\| R\right\|^2_{L^2}+\left\| \nabla R\right\|^2_{L^2}+\left\| \Delta R\right\|^2_{L^2}\right)+C  M^2 \gamma \left(\frac{3}{2}-\theta\right)^{2}e^{-4\tau}.
}
}
Using Hölder and the convexity inequality $ab\leq \frac{1}{4}a^4+\frac{3}{4}b^{\frac{4}{3}}$, a simple computation leads to
\cent{
$
\displaystyle \left\|\left|X\right|^3 R\right\|^2_{L^2} \leq \frac{3\mu^{4/3}}{4}\left\|\left|X\right|^4 R\right\|^2_{L^2}+\frac{1}{4 \mu^4} \left\| R\right\|^2_{L^2},
$
}
for all $\mu>0$.
\vspace{0.5cm}\\
Using this inequality with $\mu$ small enough, we finally obtain
\eq{\label{eqE55}\arr{l}{
\displaystyle \partial_\tau E_5+3E_5+\frac{1}{16}\left\|\left|X\right|^4 R\right\|^2_{L^2}+\left(\frac{\alpha}{2} e^{-\tau}+\frac{\alpha^2}{4} e^{-2\tau}\right)\left\|\left|X\right|^4 \Delta R\right\|^2_{L^2} \leq\displaystyle K_1\left\| R\right\|^2_{L^2}\\
\hspace{1.5cm}  + CM \gamma^{1/4} \left(\frac{3}{2}-\theta\right)^{1/2}\left(\left\| R\right\|^2_{L^2}+\left\| \nabla R\right\|^2_{L^2}+\left\| \Delta R\right\|^2_{L^2}\right)+C M^2 \gamma \left(\frac{3}{2}-\theta\right)^{2} e^{-4\tau},
}
}
where $K_1$ is a positive constant.
\begin{flushright}
$\square$
\end{flushright}
This lemma, combined with the inequality (\ref{eqE4}) enables to finish the $H^2(4)$ estimate of $R$. We define the functional
\eq{\label{E6K}
E_6 = K E_4 + E_5,
}
with $K$ some large positive constant that will be made more precise later.
\vspace{0.5cm}\\
Inequalities (\ref{eqE4}) and (\ref{eqE5}) show that one has
\cent{
$
\arr{l}{
\displaystyle\partial_\tau E_6 +2\theta E_6 + 10K\left\|\left(-\Delta\right)^{-\left(\theta-\frac{1}{4}\right)} R\right\|^2_{L^2} + \frac{K}{2}\left\|\nabla R\right\|^2_{L^2}+ \frac{K}{4}\left\|\Delta R\right\|^2_{L^2}\\
\\
\displaystyle\hspace{6cm}+\frac{1}{16}\left\|\left|X\right|^4 R\right\|^2_{L^2}+\frac{\alpha^2}{4} e^{-2\tau}\left\|\left|X\right|^4 \Delta R\right\|^2_{L^2}\leq \\
\\
\hspace{0.5cm}\displaystyle\hspace{0cm}K_1\left\| R\right\|^2_{L^2} + CM^2 \gamma^{1/4} \left(\left\| R\right\|^2_{L^2}+\left\| \nabla R\right\|^2_{L^2}+\left\| \Delta R\right\|^2_{L^2}+\left\|\left|X\right|^4 R\right\|^2_{L^2}\right)+C M^2 \gamma \left(\frac{3}{2}-\theta\right) e^{-\frac{7\tau}{2}}.
}
$
}
Interpolating again $\left\| R\right\|^2_{L^2}$ between $\left\|\left(-\Delta\right)^{-\left(\theta-\frac{1}{4}\right)} R\right\|^2_{L^2}$ and $\left\|\nabla R\right\|^2_{L^2}$ and taking $K$ and $\gamma_0$ respectively sufficiently large and small, we get
\eq{\label{eqE6}
\displaystyle\partial_\tau E_6 +2\theta E_6\leq C  M^2 \gamma \left(\frac{3}{2}-\theta\right)e^{-\frac{7\tau}{2}}.
}
\section{\label{seclim} Proof of Theorem \ref{theo1}}
\subsection{Theorem \ref{theo1} for approximate solutions}
In this section, under the condition (\ref{cond1}), we show that the solutions of (\ref{g2We}) are actually global in time and that the inequality (\ref{inetheo1}) of Theorem \ref{theo1} holds for these solutions. To get this result, we take advantage of the energy estimates that we have obtained in Section \ref{secenergy}. The following theorem is a copy of Theorem \ref{theo1} for solutions of the regularized system (\ref{g2We}).
\theo{
\label{theo2}
Let $\theta$ be a fixed positive constant such that $0<\theta<\frac{3}{2}$, $\varepsilon$ be a positive constant and $W_0 \in \mathbb H^2(4)$. There exist three positive constants $\gamma_0=\gamma_0(\alpha)$, $\varepsilon=\varepsilon_0(\alpha)$ and $T = T_0(\alpha,\theta)$ such that if $T\geq T_0$, $\varepsilon \leq \varepsilon_0$ and there exists a positive constant $\gamma\leq \gamma_0$ such that $W_0 \in \mathbb H^2(4)$ satisfies the condition
\eq{\label{cond3}
\left\|W_0\right\|^2_{L^2(4)}+ \left\|\nabla W_0\right\|^2_{L^2}+\alpha e^{-\tau_0} \left\|\Delta W_0\right\|^2_{L^2} + \alpha^2 e^{-2\tau_0}\left\|\left|X\right|^4 \Delta W_0\right\|^2_{L^2} \leq \gamma \left(\frac{3}{2}-\theta\right)^2,
}
where $\tau_0 = \log(T)$,
\vspace{0.5cm}\\
then there exist a unique solution $W_\varepsilon \in C^1\left(\left(\tau_0,+\infty\right), \mathbb H^1(4)\right)\cap C^0\left(\left(\tau_0,+\infty\right), \mathbb H^3(4)\right)$ to the system (\ref{g2We}) and a positive constant $C=C(\alpha, \tau_0)$ such that, for all $\tau \geq \tau_0$,
\eq{\label{ineapprox}
\left\|\left(Id-\alpha e^{-\tau}\Delta\right)\left(W_\varepsilon(\tau)-e^{-\tau}\sum^3_{i=1} b_i f_i\right)\right\|_{L^2(4)} \leq C \gamma \left(\frac{3}{2}-\theta\right) e^{-\theta \tau},
}
where $\displaystyle b_i = \int_{\R^3} p_i(X) . W_0( X)dX$.}
In order to prove this theorem, we use the energy estimates that we established in the section \ref{secenergy}. To obtain the inequality (\ref{ineapprox}), we need the energy functional $E_6$ to be equivalent to the $H^2(4)$-norm of $R_\varepsilon$. If we take $K$ large enough in the definition (\ref{E6K}) of $E_6$, then the next lemma holds.
\lem{\label{lemeq}
Let $R_\varepsilon \in C^1\left(\left(\tau_0,+\infty\right), \mathbb H^1(4)\right)\cap C^0\left(\left(\tau_0,+\infty\right), \mathbb H^3(4)\right)$ and $E_6$ be the energy functional defined by (\ref{E6K}). There exists $K_0$ such that, if $K\geq K_0$, then there exists a positive constant $C$ such that
\eq{\label{equival1}
E_6 (\tau) \leq C\left(\left\|R_\varepsilon\right\|^2_{L^2(4)}+\left\|\nabla R_\varepsilon \right\|^2_{L^2}+\alpha e^{-\tau}\left\|\Delta R_\varepsilon\right\|^2_{L^2} +\alpha^2 e^{-2\tau} \left\|\left|X\right|^4\Delta R_\varepsilon\right\|^2_{L^2(4)}\right),
}
\eq{\label{equival2}
C \left(\left\|R_\varepsilon\right\|^2_{L^2(4)}+\left\|\nabla R_\varepsilon \right\|^2_{L^2}+\alpha e^{-\tau}\left\|\Delta R_\varepsilon\right\|^2_{L^2} +\alpha^2 e^{-2\tau} \left\|\left|X\right|^4\Delta R_\varepsilon\right\|^2_{L^2(4)}\right) \leq E_6(\tau).
}
}
\textbf{Proof: }The inequalities (\ref{equival1}) and (\ref{equival2}) come directly from the definition of $E_6$ and the interpolation inequality (\ref{interp}).
\begin{flushright}
$\square$
\end{flushright}
\paragraph{Proof of theorem \ref{theo2}:}
Let $\theta$ be a fixed constant such that $0 <\theta<\frac{3}{2}$ and\\ $W_\varepsilon \in C^1\left(\left(\tau_0,+\infty\right), \mathbb H^1(4)\right)\cap C^0\left(\left(\tau_0,+\infty\right), \mathbb H^3(4)\right)$ be the solution of the system (\ref{g2We}) given by Theorem \ref{theoapprox}. Let $T$ and $K$ be sufficiently large so that they satisfy the conditions of the lemmas \ref{lemE0}, \ref{lemE1}, \ref{lemE3} and \ref{lemE5} and assume that the initial data $W_0$ satisfy the condition (\ref{cond1}) for some $\gamma>0$ which will be made more precise later. We decompose $W_\varepsilon$ such that
\cent{
$
W_\varepsilon = e^{-\tau} \Omega_\infty + R_\varepsilon,
$
}
where $\Omega_\infty= \displaystyle \sum^3_{i=1} b_i f_i$, $\displaystyle b_i = \int_{\R^3} p_i(X).W_0(X) dX$ and $\left\{f_1,f_2,f_3\right\}$ is the basis of the eigenspace of $\mathcal L$ associated to the eigenvalue $-1$, given by (\ref{fi}).
\vspace{0.5cm}\\
Let $M$ be a positive constant such that $M>2$ that will be made more precise later and $\tau^*_\varepsilon \in \left[\tau_0,\tau_\varepsilon\right]$ be the biggest positive time such that the inequality (\ref{cond2W}) holds. We take $\gamma$ and $\varepsilon$ sufficiently small so that the lemmas \ref{lemE0}, \ref{lemE1}, \ref{lemE3} and \ref{lemE5} hold. According to the inequality (\ref{eqE6}), one has, for all $\tau \in \left[\tau_0,\tau^*_\varepsilon\right)$,
\eq{
\displaystyle \partial_\tau \left(E_6(\tau) e^{2\theta\tau}\right) \leq C M^2 \gamma\left(\frac{3}{2}-\theta\right)  e^{-\left(\frac{7}{2}-2\theta\right)\tau }
}
Integrating in time the previous inequality between $\tau_0$ and $\tau \in \left[\tau_0, \tau^*_\varepsilon\right)$, we get
\eq{\label{eqpreuve1}
\displaystyle E_6(\tau) \leq E_6(\tau_0) e^{ -2\theta\left(\tau -\tau_0\right)} + C M^2 \gamma \left(\frac{3}{2}-\theta\right) e^{-\frac{7\tau_0}{2}}\left(e^{-2\theta \left(\tau-\tau_0\right)}-e^{-\frac{7}{2} \left(\tau-\tau_0\right)}\right).
}
Arguing like in the proof of Lemma \ref{conditions} and using the inequality (\ref{equival1}), we can show that
\cent{
$
\displaystyle E_6(\tau_0) \leq C_1 \gamma \left(\frac{3}{2}-\theta\right)^2,
$
}
which implies
\eq{
E_6 (\tau) \leq C\gamma\left(\frac{3}{2}-\theta\right)^2 + C M^2 \gamma \left(\frac{3}{2}-\theta\right) e^{-\frac{7\tau_0}{2}}.
}
According to the inequalities (\ref{equival2}) and (\ref{cond2}), one has, for all $\tau \in \left[\tau_0,\tau^*_\varepsilon\right)$,
\cent{$
\arr{l}{
\displaystyle \left|b\right|^2+\left\|R_\varepsilon\right\|^2_{L^2(4)}+\left\|\nabla R_\varepsilon\right\|^2_{L^2} + \alpha e^{-\tau}\left\|\Delta R_\varepsilon\right\|^2_{L^2} +\alpha^2 e^{-2\tau} \left\|\left|X\right|^4\Delta R_\varepsilon \right\|^2 \leq \\
\hspace{8cm}\displaystyle C\gamma\left(\frac{3}{2}-\theta\right)^2 + C M^2 \gamma \left(\frac{3}{2}-\theta\right) e^{-\frac{7\tau_0}{2}}.
}
$}
Recalling that $\displaystyle W_\varepsilon = \sum^3_{i=1} b_i f_i + R_\varepsilon$, we get
\cent{$
\arr{l}{
\displaystyle \left\|W_\varepsilon\right\|^2_{L^2(4)}+\left\|\nabla W_\varepsilon\right\|^2_{L^2} + \alpha e^{-\tau}\left\|\Delta W_\varepsilon\right\|^2_{L^2} +\alpha^2 e^{-2\tau} \left\|\left|X\right|^4\Delta W_\varepsilon \right\|^2 \leq \\
\hspace{8cm} \displaystyle C_1\gamma\left(\frac{3}{2}-\theta\right)^2 + C_2 M^2 \gamma \left(\frac{3}{2}-\theta\right) e^{-\frac{7\tau_0}{2}},
}
$}
where $C_1$ and $C_2$ are two positive constants.
\vspace{0.5cm}\\
We take $M$ sufficiently large so that $C_1\leq \frac{M}{4}$ and $\tau_0 = \ln(T)$ sufficiently large so that $\displaystyle C_2 M^2 e^{-\frac{7\tau_0}{2}} \leq \frac{M \left(\frac{3}{2}-\theta\right)}{4}$, we obtain, for all $\tau \in \left[\tau_0,\tau^*_\varepsilon\right)$,
\eq{\label{eqpreuve2}
\displaystyle \left\|W_\varepsilon\right\|^2_{L^2(4)}+\left\|\nabla W_\varepsilon\right\|^2_{L^2} + \alpha e^{-\tau}\left\|\Delta W_\varepsilon\right\|^2_{L^2} +\alpha^2 e^{-2\tau} \left\|\left|X\right|^4\Delta W_\varepsilon \right\|^2 \leq  \frac{M\gamma\left(\frac{3}{2}-\theta\right)^2}{2}.
}
In particular, the inequality (\ref{eqpreuve2}) shows that $\tau^*_\varepsilon=\tau_\varepsilon$. Furthermore, letting $\tau$ tend to $\tau_\varepsilon$, we see that if $\tau_\varepsilon$ is finite, then the $H^1(4)$ norm of $W_\varepsilon$ stay bounded on $\left[\tau,\tau_\varepsilon\right)$. According to the proof of Theorem \ref{theoapprox}, it implies in particular that one can extend the interval of definition of $W_\varepsilon$ over $\tau_\varepsilon$. Consequently, we have necessarily $\tau_\varepsilon = +\infty$. In addition, going back to the inequality (\ref{eqpreuve1}) and applying the inequality (\ref{equival2}) of Lemma \ref{lemeq}, we see that the inequality (\ref{ineapprox}) holds.
\begin{flushright}
$\square$
\end{flushright}
\subsection{\label{subseclim} Existence of weak solutions in $\mathbb H^2(4)$}
In this section, we show that there exists a weak solution to the system (\ref{g2W}) belonging to the space $C^0\left(\left[\tau_0,+\infty\right),\mathbb H^2(4)\right)$. To this end, we show that, when $\varepsilon$ tends to $0$, $W_\varepsilon$ tends to a divergence free vector $W$ which satsifies (\ref{g2W}) in a weak sense. Let $\left(\varepsilon_n\right)_{n\in \mathbb N}$ be a sequence of positive terms which tends to $0$. Let $W_{\varepsilon_n} \in C^1\left(\left(\tau_0,+\infty\right), \mathbb H^1(4)\right) \cap C^0\left(\left(\tau_0,+\infty\right), \mathbb H^3(4)\right)$ be the global solution of (\ref{g2We}) given by Theorem \ref{theo2}, with initial data $W_0$. Let $\mathcal O$ be a bounded open set of $\R^3$. For $s \in \R^+$, $H^s(\mathcal O)$ denotes the restriction of the Sobolev space $H^s(\R^3)$ on $\mathcal O$. For $s \geq 1$, we define also the space 
\cent{
$
\displaystyle H^s_0(\mathcal O)=\left\{u \in H^s(\mathcal O):u_{\left|\partial \mathcal O\right.=0}\right\}.
$
}
Let $\tau_1$ be a fixed positive time such that $\tau_1 >\tau_0$. Due to the boundedness property of $W_{\varepsilon_n}$  in $L^\infty\left(\left[\tau_0,\tau_1\right], \mathbb H^2(4)\right)$ uniformly with respect to $n$, there exist $W \in L^\infty\left(\left[\tau_0,\tau_1\right], \mathbb H^2(4)\right)$ and a subsequence of $\varepsilon_n$ (that we still note $\varepsilon_n$) such that 
\eq{\label{idlim1}
W_{\varepsilon_n} \rightharpoonup W\quad \hbox{ weak* in} \quad L^\infty\left(\left[\tau_0,\tau_1\right],H^2(\mathcal O)^3\right).
}
Since $W_{\varepsilon_n}$ is bounded in $L^\infty\left(\left[\tau_0,\tau_1\right], \mathbb H^2(4)\right)$, applying the operator $\left(I-\alpha e^{-\tau} \Delta\right)^{-1}$ to the first equality of (\ref{g2We}), it is quite easy to see that $\partial_\tau W_{\varepsilon_n}$ is bounded in $L^\infty\left(\left[\tau_0,T\right],L^2(\mathcal O)^3\right)$ uniformly with respect to $n$. Consequently, $W_{\varepsilon_n}$ is equicontinuous in time on $L^2(\mathcal O)^3$. Indeed, given $\sigma_1$ and $\sigma_2$ belonging to $\left[\tau_0,\tau_1\right]$, one has
\aligne{
\left\|W_{\varepsilon_n}(\sigma_1) -W_{\varepsilon_n}(\sigma_2)\right\|_{L^2(\mathcal O)} &= \left\|\int^{\sigma_1}_{\sigma_2}\partial_\tau W_{\varepsilon_n}(s)ds \right\|_{L^2(\mathcal O)}\\
& \leq \left|\int^{\sigma_1}_{\sigma_2}\left\|\partial_\tau W_{\varepsilon_n}(s) \right\|_{L^2(\mathcal O)}ds\right|\\
&\leq \left|\sigma_1-\sigma_2\right|\max\limits_{s\in \left[\tau_0,T\right]}\left\|\partial_\tau W_{\varepsilon_n}(s) \right\|_{L^2(\mathcal O)}.
}
Besides, for all $\tau \in \left[\tau_0,\tau_1\right]$, the set $\bigcup\limits_{n\in \mathbb N} W_{\varepsilon_n}(\tau)$ is bounded in $H^2(\mathcal O)^3$ and thus compact in $L^2(\mathcal O)^3$. Applying the classical Arzela-Ascoli theorem, we conclude that
\cent{
$W_{\varepsilon_n} \longrightarrow W$ strongly in $C^0\left(\left[\tau_0,\tau_1\right],L^2(\mathcal O)^3\right)$.
}
A classical interpolation inequality between $L^2$ and $H^2$ yields, for all $s<2$,
\eq{\label{idlim2}
W_{\varepsilon_n} \longrightarrow W \quad \hbox{strongly in} \quad C^0\left(\left[\tau_0,\tau_1\right],H^s(\mathcal O)^3\right).
}
The two identities (\ref{idlim1}) and (\ref{idlim2}) are sufficient to pass to the limit in the weak formulation of the system (\ref{g2We}) and to show that $W$ is a weak solution of the system (\ref{g2W}). More precisely, for every $\varphi \in C^1\left(\left[\tau_0,\tau_1\right],H^1_0(\mathcal O)^3\right)$ such that $\dv \varphi =0$, one has, for all $\tau \in \left[\tau_0,\tau_1\right]$,
\eq{\label{weakg2}
\arr{l}{
\displaystyle \int_{\mathcal O} \left(W(\tau)-\alpha e^{-\tau}\Delta W(\tau) \right).\varphi(\tau) dX+\int^\tau_{\tau_0}\int_{\mathcal O} \mathcal L\Big(W(\sigma)\Big) .\varphi (\sigma) dX d\sigma \\
\\
\displaystyle \hspace{5cm}+\int^\tau_{\tau_0}\int_{\mathcal O} \left(W(\sigma)-\alpha e^{-\sigma}\Delta W(\sigma)\right)\wedge U(\sigma).\curl \varphi (\sigma) dX d\sigma\\
\\
\displaystyle = \int_{\mathcal O} \left(W_0-\alpha e^{-\tau_0}\Delta W_0 \right).\varphi(\tau_0) dX +\int^\tau_{\tau_0}\int_{\mathcal O} \left(W(\sigma)-\alpha e^{-\sigma}\Delta W(\sigma) \right).\partial_\tau \varphi(\sigma) dX d\sigma\\
\\
\displaystyle \hspace{2cm}+\int^\tau_{\tau_0}\int_{\mathcal O} \frac{3\alpha}{2}e^{-\sigma} \Delta W(\sigma).\varphi(\sigma) dX d\sigma +\int^\tau_{\tau_0}\int_{\mathcal O} \frac{\alpha}{2} e^{-\sigma} \Delta W(\sigma)\left(X.\nabla \varphi(\sigma)\right) dX d\sigma .
}
}
We just show that the non-linear term converges, using (\ref{idlim1}) and (\ref{idlim2}). The other ones are nearly obvious. We have
\eq{
\arr{l}{\displaystyle \int^\tau_{\tau_0}\int_{\mathcal O} \left(W_{\varepsilon_n}(\sigma)-\alpha e^{-\sigma}\Delta W_{\varepsilon_n}(\sigma)\right)\wedge U_{\varepsilon_n}(\sigma).\curl \varphi (\sigma) dX d\sigma =\\
\displaystyle \hspace{3cm} \int^\tau_{\tau_0}\int_{\mathcal O} \left(W(\sigma)-\alpha e^{-\sigma}\Delta W(\sigma)\right)\wedge U (\sigma).\curl \varphi (\sigma) dX d\sigma+R_n +S_n,
}
}
where
\cent{
$
\arr{l}{
R_n = \displaystyle \int^\tau_{\tau_0}\int_{\mathcal O} \left(W_{\varepsilon_n}(\sigma)-\alpha e^{-\sigma}\Delta W_{\varepsilon_n}(\sigma)\right)\wedge \left(U(\sigma)-U_{\varepsilon_n}(\sigma)\right).\curl \varphi (\sigma) dX d\sigma,\\
\\
S_n = \displaystyle \int^\tau_{\tau_0}\int_{\mathcal O}\left(W(\sigma)-W_{\varepsilon_n}(\sigma)-\alpha e^{-\sigma}\left(\Delta W(\sigma)-\Delta W_{\varepsilon_n}(\sigma)\right)\right) \wedge U(\sigma).\curl \varphi (\sigma) dX d\sigma.
}
$
}
Due to Hölder inequalities, the boundedness property of $W_{\varepsilon_n}$ in $H^2(\mathcal O)^3$ and the inequality (\ref{biotsb}), we have
\aligne{
R_n &\leq C \int^\tau_{\tau_0} \left\|U(\sigma)-U_{\varepsilon_n}(\sigma)\right\|_{L^\infty(\mathcal O)}\left\|\nabla \varphi(\sigma)\right\|_{L^2(\mathcal O)}d\sigma\\
&\leq C\int^\tau_{\tau_0} \left\|W(\sigma)-W_{\varepsilon_n}(\sigma)\right\|^{1/2}_{L^2(\mathcal O)}\left\|W(\sigma)-W_{\varepsilon_n}(\sigma)\right\|^{1/2}_{L^6(\mathcal O)}\left\|\nabla \varphi(\sigma)\right\|_{L^2(\mathcal O)}d\sigma\\
&\leq C\left(T-\tau_0\right)\max\limits_{\sigma \in \left[\tau_0,T\right]} \left\|W(\sigma)-W_{\varepsilon_n}(\sigma)\right\|_{H^1(\mathcal O)}\max\limits_{\sigma \in \left[\tau_0,T\right]}\left\|\nabla \varphi(\sigma)\right\|_{L^2(\mathcal O)}.
}
Thus, the identity (\ref{idlim2}) implies that $R_n \rightarrow 0$ when $n\rightarrow +\infty$.
\vspace{0.5cm}\\
Because of the identity (\ref{idlim1}), it is clear that we have also $S_n \rightarrow 0$ when $n\rightarrow +\infty$. Thus, we have shown that, for all $\tau \in \left[\tau_0,\tau_1\right]$,
\eq{
\arr{l}{\displaystyle \lim\limits_{n\rightarrow +\infty} \int^\tau_{\tau_0}\int_{\mathcal O} \left(W_{\varepsilon_n}(\sigma)-\alpha e^{-\sigma}\Delta W_{\varepsilon_n}(\sigma)\right)\wedge U_{\varepsilon_n}(\sigma).\curl \varphi (\sigma) dX d\sigma =\\
\\
\displaystyle \hspace{5cm} \int^\tau_{\tau_0}\int_{\mathcal O} \left(W(\sigma)-\alpha e^{-\sigma}\Delta W(\sigma)\right)\wedge U(\sigma).\curl \varphi (\sigma) dX d\sigma.
}
}
Furthermore, since $W_{\varepsilon_n}(\tau)$ converge weakly to $W(\tau)$ in $\mathbb H^2(4)$, from the inequality (\ref{ineapprox}), we get
\eq{
\left\|\left(I-\alpha e^{-\tau} \Delta\right)\left(W (\tau)- e^{-\tau}\sum^3_{i=1} b_i f_i\right)\right\|_{L^2(4)} \leq C \gamma \left(\frac{3}{2}-\theta\right) e^{-\theta\tau},
}
for all $\tau \in \left[\tau_0,+\infty\right)$.
\subsubsection{Uniqueness}
It remains to show that the solutions of (\ref{g2w}) are unique in the space $C^0\left(\left[0,+\infty\right),\mathbb H^2(4)\right)$. To show this fact, it suffices to show that the divergence free vector field $u$ obtained from a solution $w$ of (\ref{g2w}) through the Biot-Savart law is unique. Since $w$ belongs to $C^0\left(\left[0,+\infty\right),\mathbb H^2(4)\right)$, the inequality (\ref{biotsa}) with $q=2$ and $p=\frac{6}{5}$ and the inequality (\ref{biotsc}) with $p=2$ of the lemma \ref{biots} imply directly that $u \in C^0\left(\left[0,+\infty\right),H^3(\R^3)^3\right)$. Furthermore, $u$ satisfies the equations of motion of second grade fluids (\ref{g2u}). The uniqueness of the $H^3-$solutions of (\ref{g2u}) has been shown in \cite{cioranescugirault97} for the case of a bounded open set of $\R^3$ with Dirichlet boundary conditions. In our case, we can apply the computations of the proof of \cite[Theorem 2]{cioranescuelhacene84}, which imply the uniqueness of the solutions of (\ref{g2u}) with initial data in $H^3(\R^3)^3$.

\end{document}